\documentclass[moor]{informs4post}
\usepackage{eqndefns-left} % For checking the display equation width and equation environment definitions %
\RequirePackage{tgtermes}
\RequirePackage{newtxtext}
\RequirePackage{newtxmath}
\RequirePackage{bm}
\RequirePackage{endnotes}

\SingleSpacedXI
\usepackage[sort&compress]{natbib}
 \bibpunct[, ]{[}{]}{,}{n}{}{,}%
 \def\bibfont{\small}%
\EquationsNumberedBySection % (1.1), (1.2), ...

\TheoremsNumberedBySection  % (Theorem 1.1, Lemma 1.1, Theorem 1.2)
\ECRepeatTheorems  %  

\MANUSCRIPTNO{MOR-2024-0853}

\usepackage[english]{babel}
\usepackage[utf8]{inputenc} % allow utf-8 input
\usepackage[T1]{fontenc}    % use 8-bit T1 fonts
\usepackage{url}            % simple URL typesetting
\usepackage{booktabs}       % professional-quality tables
\usepackage{amsfonts}       % blackboard math symbols
\usepackage{nicefrac}       % compact symbols for 1/2, etc.
\usepackage{microtype}      % microtypography
\usepackage[margin=1in]{geometry}
\usepackage[ruled,linesnumbered]{algorithm2e}
\usepackage{amssymb}
\usepackage{amsfonts}
\usepackage{amsmath}
\usepackage{subcaption}

\usepackage{pdfcomment}
\usepackage{array}
\usepackage{enumitem}
\usepackage{graphicx} 

{\theoremstyle{THkey}\newtheorem{condition}{Condition}}
{\theoremstyle{THkey}}

\newcommand{\ignore}[1]{}

\renewcommand{\Pr}{{\rm Pr}}

\newcommand{\E}{{\rm E}}

\newcommand{\eps}{\varepsilon}
\newcommand{\R}{\mathbb R}

\newcommand{\simplephi}{\Phi}

\newcommand{\calS}{{\cal S}}

\newcommand{\calX}{{\cal X}}
\newcommand{\calY}{{\cal Y}}

\newcommand{\x}{{\boldsymbol x}}

\begin{document}
%%%%%%%%%%%%%%%%

%\setcounter{page}{1}
%\DOI{10.1287/xxxx.0000.0000}%

\RUNAUTHOR{Xiong}

% \RUNTITLE{Predictive Maintenance in Manufacturing}
\RUNTITLE{High-Probability Polynomial-Time Complexity of Restarted PDHG for Linear Programming}

% \TITLE{Optimal Resource Allocation in Humanitarian Logistics: A Stochastic Programming Approach}
\TITLE{High-Probability Polynomial-Time Complexity of Restarted PDHG for Linear Programming}

\ARTICLEAUTHORS{%
%\AUTHOR{John Doe,\textsuperscript{a} Jane Smith,\textsuperscript{b}}
%\AFF{\textsuperscript{a}Department of Industrial Engineering, University of XYZ, \EMAIL{john.doe@xyz.edu; \textsuperscript{b}Department of Computer Science, University of ABC, \EMAIL{jane.smith@abc.edu}} 
\AUTHOR{Zikai Xiong}
\AFF{Department of Industrial Engineering and Management Sciences, Northwestern University, 2145 Sheridan Road, Evanston, IL 60208, USA. \href{mailto:zikai.xiong@northwestern.edu}{zikai.xiong@northwestern.edu}.}
} % end of the block

\ABSTRACT{%
The restarted primal-dual hybrid gradient method (rPDHG) is a first-order method recently known for its computational effectiveness in solving linear programming (LP) problems.  Despite its impressive practical performance, the theoretical iteration bounds for rPDHG can be exponentially poor. To investigate this gap from a probabilistic perspective, we show that rPDHG achieves polynomial-time complexity in a high-probability sense, under assumptions on the probability distribution from which the data instance is generated. We consider both Gaussian and more general sub-Gaussian distribution models. For standard-form LP instances with $m$ constraints and $n$ variables, our bounds take a particularly simple form when $m$ is not too close to $n$: rPDHG iterates settle on the optimal basis in $\widetilde{O}\left(\tfrac{n^{2.5}m^{0.5}}{\delta}\right)$ iterations, followed by $O\left(\frac{n^{0.5}m^{0.5}}{\delta}\ln\big(\tfrac{1}{\varepsilon}\big)\right)$ iterations to compute an $\varepsilon	$-optimal solution,  with probability at least $1-\delta$ for $\delta$ that is not exponentially small. The Stage-I bound further improves to $\widetilde{O}\left(\frac{n^{2.5}}{\delta}\right)$ in the Gaussian distribution model.  Experimental results confirm the tail behavior and the polynomial-time dependence on problem dimensions of the iteration counts. As an application of our probabilistic analysis, we explore how the disparity among the components of the optimal solution bears on the performance of rPDHG, and we provide guidelines for generating challenging LP instances. 
}%

%Supplemental Material:
%Data Ethics & Reproducibility Note:

%\KEYWORDS{Stochastic programming, Decision support,Uncertainty, Disaster response, Optimization}

\renewcommand{\MSCCLASSname}{{\bf MSC2020 subject classification\/}{\kern0.7pt}:\enskip}
\MSCCLASS{Primary: 90C05; Secondary: 90C60}
% \ORMSCLASS{Primary: Programming/Linear; Secondary: Analysis of Algorithms/Probabilistic} 
% \HISTORY{Received December 27, 2024; revised February 14, 2026, and August 2, 2026; accepted September 8, 2026.}

\maketitle
%%%%%%%%%%%%%%%%%%%%%%%%%%%%%%%%%%%%%%%%%%%%%%%%%%%%%%%%%%%%%%%%%%%%%%

\section{Introduction} 
 
Linear programming (LP) has been one of the most fundamental areas of optimization since the 1950s, with applications spanning many domains. For over 70 years, researchers and practitioners have extensively studied and refined LP solution methods, primarily through two major algorithmic frameworks: the simplex method, introduced by Dantzig in 1947, and the interior-point method, developed by Karmarkar in 1984. 
The profound impact of these methods on both academia and industry cannot be overstated.

While the two classic methods (simplex and interior-point methods) are successful in solving many LP problems, both of them rely on frequent matrix factorizations, whose computational cost grows superlinearly with problem dimensions. Moreover, these matrix factorizations cannot efficiently utilize modern computational architectures, especially graphics processing units (GPUs). These limitations have motivated the development of the restarted primal-dual hybrid gradient method (rPDHG), a first-order method that eliminates the need for matrix factorizations. By utilizing primarily matrix-vector products for gradient computations, rPDHG achieves better scalability by exploiting problem sparsity and leveraging parallel architectures. It directly addresses the saddle-point formulation (see \citet{applegate2023faster}),  and also solves conic linear programs (see \citet{xiong2024role}) and convex quadratic programs (see \citet{lu2023practicalqpnew,huang2024restarted}).

The rPDHG method has demonstrated performance comparable to, and sometimes exceeding, that of classical simplex and interior-point methods on many LP instances. Its implementations span both CPUs (PDLP by \citet{applegate2021practical}) and GPUs (cuPDLP and cuPDLP-C by \citet{lu2023cupdlpnew} and \citet{lu2023cupdlp-c}). Recognizing its effectiveness, leading commercial solvers including COPT 7.1 (see \citet{coptgithub}), Xpress 9.4 (see \citet{xpressnwes}), and Gurobi 13.0 (see \citet{gurobinews}) have integrated rPDHG as a third foundational LP algorithm alongside simplex and interior-point methods. The method has also been adopted by other solvers such as Google OR-Tools (see \citet{applegate2021practical}), HiGHS (see \citet{highsrelease}), and NVIDIA cuOpt (see \citet{nvdianews}).

A persistent challenge in LP algorithm development has been the ``embarrassing gap'' between the practical experience of these LP algorithms' efficiency and their worst-case iteration bounds. The efficiency is typically measured by the iteration count (to solve the problem) as a function of problem dimensions. For simplex methods,  empirical evidence suggests it is usually a low-degree polynomial or even a linear function of the dimensions (see, for example, \citet{shamir1987efficiency}), but worst-case analyses reveal this function could potentially be exponentially large (see, for example, \citet{klee1972good}).  For interior-point methods, worst-case iteration bounds are usually a polynomial of the dimension times the ``bit-length'' $L$ of the problem, but the practical performance is usually much better than the worst-case analyses would suggest. This gap between theory and practice has motivated extensive research into average-case complexity analysis for both classic LP methods, which we review shortly after.

A similar gap exists between the theoretical bound and empirical performance for rPDHG. 
While the method demonstrates strong practical performance on most LP problems (see, e.g., \citet{applegate2021practical,lu2023cupdlpnew}), its theoretical worst-case iteration bounds can be exponentially poor. Recent work has made significant progress in understanding rPDHG's convergence behavior. \citet{applegate2023faster} establish the first linear convergence rate of rPDHG using the global Hoffman constant of the matrix of the KKT system that defines the optimal solutions.  \citet{xiong2023computational,xiong2023relation} provide a tighter computational guarantee for rPDHG using two natural properties of the LP problem: LP sharpness and the ``limiting error ratio.''   Furthermore, \citet{xiong2024accessiblenew} gives a closed-form iteration bound for LPs with unique optima and demonstrates that rPDHG has a two-stage performance. In Stage I, the iterates settle on the optimal basis (and thus identify this basis). Stage II then exhibits faster local convergence to the optimal solution. However, all of the existing computational guarantees at least linearly depend on certain condition measures that can be exponentially large in the worst case, which leaves unexplained the strong practical performance of rPDHG observed in many applications.
 
To help address this gap between theory and practice from a probabilistic perspective, this paper primarily aims at answering the following question:
\vspace{8pt}
\begin{enumerate}
	\item[Q1.] Can we prove that rPDHG is efficient with high probability under a probabilistic input model?
\end{enumerate}
\vspace{8pt}
Here, ``efficient'' means having polynomial-time complexity. An affirmative answer to this question would provide a theoretical explanation for rPDHG's strong performance under the probabilistic model and, more importantly, offer insights into its behavior beyond the worst case. We approach this question through probabilistic analysis and give an affirmative answer in this paper.
 
Probabilistic analysis has successfully narrowed similar theory-practice gaps for classic LP methods. By assuming the probability distributions according to which the problem data was generated (or equivalently, the distributions on the input data), researchers can provide bounds on the expected number of iterations. For the simplex method, breakthrough results in the 1980s established various polynomial iteration bounds under various probabilistic assumptions (see, for example,  \citet{borgwardt1977b,borgwardt1978,borgwardt1982average,borgwardt1982some,smale1983average,smale1983problemnew,haimovich1983simplex,adler1985simplex,adler1986family,adler1987simplex,megiddo1986improved}).  Similarly, for interior-point methods, \citet{todd1991probabilistic} proposes a probabilistic model. Various versions of this model lead to polynomial expected iteration bounds in terms of problem dimensions, independent of the bit-length $L$ (see, for example, \citet{ye1994toward,anstreicher1993averagenew,anstreicher1999probabilistic,huang2005expected,huang2000average,ji2008probabilistic}). More recently, another probabilistic analysis approach called smoothed analysis provides a new framework that bridges worst-case and average-case analyses.  Polynomial-time complexity bounds have also been established under this framework for both simplex and interior-point methods for LP instances (see, for example, \citet{spielman2004smoothed,spielman2003smoothed,dunagan2011smoothed,dadush2018friendly}).

Probabilistic analysis may also provide new tools for revealing new insights into rPDHG. Some practical observations of rPDHG, such as the two-stage performance and sensitivity under perturbations, have been partially explained by the worst-case iteration bounds established by \citet{xiong2024accessiblenew} and others. Probabilistic analysis may provide complementary perspectives and clearer insights into these behaviors of rPDHG. 
 
Despite the above, it should be admitted that real-life LP problems often differ significantly from the probabilistic models that people assume for probabilistic analyses. For example, almost all the existing probabilistic models are based on either the Gaussian distribution or the sign-invariant distribution (see, e.g., \citet{shamir1987efficiency} and \citet{todd1991probabilistic}). These distributional assumptions rarely match the data generation processes of practical LP instances. Consequently, extending these models to broader distribution classes remains an important research direction (see the survey by \citet{shamir1987efficiency}).

Therefore, beyond the main question Q1, this paper also addresses two additional questions:
\vspace{8pt}
\begin{enumerate}[noitemsep]
	\item[Q2.] Can the probabilistic analysis be extended beyond Gaussian distributions to more general distribution families?
	\item[Q3.] What novel insights into rPDHG's performance can be gained through probabilistic analysis?
\end{enumerate}
\vspace{8pt} 

To address the above three questions, we study a probabilistic model that is built on the classic probabilistic model proposed by \citet{todd1991probabilistic}.  Our model considers standard-form LP instances where the components of the constraint matrix follow specified Gaussian or, more generally, sub-Gaussian distributions, while the right-hand side and the objective vectors are constructed based on random primal and dual solutions.
This approach builds on \citet{todd1991probabilistic}'s framework, which has been dominant in probabilistic analyses of interior-point methods (see, e.g., \citet{anstreicher1993averagenew,anstreicher1999probabilistic} and \citet{ye1994toward}).  Compared with existing probabilistic models, our work extends the distribution of the input data from Gaussian to general sub-Gaussian distributions, which include Bernoulli and bounded distributions. This extension substantially broadens the class of random LP instances covered by our analysis. Our analysis leverages recent advances in nonasymptotic random matrix theory developed over the past two decades (see, e.g., the ICM 2010 lecture by \citet{rudelson2010non}).

\subsection{Outline and contributions}
In this paper, we consider LP instances in standard form:
\begin{equation}\label{pro:primal LP}
	\begin{aligned}
		\min_{x\in \R^n}  \quad c^\top x   \quad \ \ \text{s.t.}     \  Ax = b \, ,  \ x\geq 0
	\end{aligned}
\end{equation}
where the constraint matrix $A\in\R^{m\times n}$ and $m\le n$. Thus, the problem has $m$ linear equality constraints and $n$ decision variables. We also denote by $d := n - m$ their difference. Any LP instance can be reformulated equivalently in the standard form \eqref{pro:primal LP}. The dual to the problem \eqref{pro:primal LP} can be expressed as:
\begin{equation}\label{pro:dual LP}
	\begin{aligned}
		\max_{y\in \R^m, s\in\R^n}  \quad b^\top y  \quad \ \
		\text{s.t.} \   A^\top y + s = c \, ,\ s\geq 0 \
	\end{aligned}
\end{equation}
with $s\in \R^n$ denoting the slack.  
We let $x^\star$ be any optimal solution of \eqref{pro:primal LP} and let $s^\star$ be any optimal dual slack of \eqref{pro:dual LP}. The rest of the paper is organized as follows.

Section \ref{sec:preliminary} introduces the rPDHG algorithm for solving linear programs.

Section \ref{sec:high-probability} presents our main result: the high-probability polynomial-time complexity of rPDHG. We begin by introducing our probabilistic model with sub-Gaussian input data in Section \ref{subsec:randomLPmodel}. Building on rPDHG's two-stage behavior, in which Stage I settles on the optimal basis and Stage II achieves fast local convergence, Section \ref{subsec:probabilistic_analysis_subgaussian} presents high-probability iteration bounds for both stages. Section \ref{subsec:probabilistic_analysis_gaussian} then presents improved bounds under Gaussian input data. Table \ref{tab:complexity} contains a preview of these results in the case where $m$ is not too close to $n$. Here $\eps$ denotes the target error tolerance of the desired solution. These results establish high-probability polynomial-time complexity guarantees for rPDHG (addressing Q1),  and extend the probabilistic LP complexity analysis to sub-Gaussian input data (addressing Q2).

\begin{table}[htbp]
\TABLE
{\normalsize High-probability iteration bounds for computing an $\eps$-optimal solution using rPDHG (in the case where $m$ is not too close to $n$).\label{tab:complexity}}{\renewcommand{\arraystretch}{1.5}
\begin{tabular}{@{\hspace{20pt}}c@{\hspace{20pt}}c@{\hspace{20pt}}c@{\hspace{20pt}}c@{\hspace{20pt}}}
\hline \vspace{0.5ex}
\rule{0pt}{3ex}    & Distribution & Stage I & Stage II \\ 
\hline \vspace{1.5ex}
\rule{0pt}{5ex}\begin{tabular}{@{}c@{}} \renewcommand{\arraystretch}{0.3}Section \ref{subsec:probabilistic_analysis_subgaussian}  \vspace{-8pt}\\ (Theorem \ref{thm:complexity-random}) \end{tabular}
& sub-Gaussian 
& $\displaystyle \widetilde{O}\left(\frac{n^{2.5}m^{0.5}}{\delta}\right)$ 
& $\displaystyle O\left(\frac{n^{0.5}m^{0.5}}{\delta}\cdot\ln\left(\frac{1}{\eps}\right)\right)$ \\ 
\hline   \vspace{1.5ex}
\rule{0pt}{5ex}\begin{tabular}{@{}c@{}} \renewcommand{\arraystretch}{0.3}Section \ref{subsec:probabilistic_analysis_gaussian}   \vspace{-8pt}\\ (Theorem \ref{thm:complexity-gaussian}) \end{tabular}
& Gaussian 
& $\displaystyle \widetilde{O}\left(\frac{n^{2.5}}{\delta}\right)$ 
& $\displaystyle O\left(\frac{n^{0.5}m^{0.5}}{\delta}\cdot\ln\left(\frac{1}{\eps}\right)\right)$ \\
\hline   \vspace{1.5ex}
\rule{0pt}{5ex}\begin{tabular}{@{}c@{}} \renewcommand{\arraystretch}{0.3}Section \ref{sec:effect_of_optimal_solution}  \vspace{-8pt}\\ (Theorem \ref{thm:complexity-conditioned}) \end{tabular}
& sub-Gaussian
& $\displaystyle \widetilde{O}\left(\frac{n^{1.5}m^{0.5}\phi}{\delta}\right)$ 
& $\displaystyle O\left(\frac{n^{0.5}m^{0.5}}{\delta}\cdot\ln\left(\frac{1}{\eps}\right)\right)$ \\
\hline
\end{tabular}}{ \small
		\textit{Note.} These iteration bounds hold with probability at least $1-\delta$ when $\delta$ is not exponentially small in terms of $m$ and $n-m$. In the third row, the probabilities are conditional on the primal and dual problems having unique optimal solutions, and $\phi:=\frac{\frac{1}{n}\cdot\sum_{i=1}^n \left(x_i^\star + s_i^\star\right)}{\min_{1 \le i \le n} (x_i^\star + s_i^\star)}$ denotes the disparity ratio among the components of these unique optima. Here $O(\cdot)$ denotes bounds up to distribution-related constant factors, and $\widetilde{O}(\cdot)$ hides additional logarithmic factors.}
\end{table}

Sections \ref{sec:prove_thm_random} and \ref{sec:refined_analysis_gaussian} present the proofs of our results in Sections \ref{subsec:probabilistic_analysis_subgaussian} and \ref{subsec:probabilistic_analysis_gaussian}, respectively.

Section \ref{sec:experiments} presents the experimental results. They confirm the tail behavior of the iteration counts for both stages (the linear dependence on $\frac{1}{\delta}$), and also validate the polynomial dependence on problem dimensions in the high-probability iteration counts.  

Section \ref{sec:effect_of_optimal_solution} provides new insights using probabilistic analysis (addressing Q3). We investigate how the disparity ratio $\phi:=\frac{\frac{1}{n}\cdot\sum_{i=1}^n \left(x_i^\star + s_i^\star\right)}{\min_{1 \le i \le n} (x_i^\star + s_i^\star)}$ among the optimal solution components of $x^\star,s^\star$ influences rPDHG's performance, by deriving a high-probability iteration bound conditioned on unique primal and dual optima that grows linearly in this disparity ratio $\phi$. This result yields practical insights for generating challenging test LP instances, which we validate experimentally.

\subsection{Related work}\label{subsec:related-work}
In addition to the above worst-case analysis of rPDHG (\citet{applegate2023faster,xiong2023computational,xiong2024accessiblenew}),  
\citet{hinder2023worstnew} proves that rPDHG has a polynomial iteration bound for total unimodular LPs. \citet{xiong2024role} provide computational guarantees of rPDHG for general conic linear programs based on geometric measures of the primal-dual (sub)level set. \citet{lu2024geometrynew,lu2024restarted} study the vanilla PDHG and the Halpern restarted PDHG using trajectory-based analysis, and demonstrate the two-stage convergence characterized by Hoffman constants of a reduced linear system defined by the limiting iterate. 
These deterministic worst-case analyses apply broadly, but these iteration bounds are parameterized by instance-specific constants that can be exponentially large in the worst case. The focus of this paper is complementary: we provide a probabilistic explanation of rPDHG's favorable behavior on a broad class of random instances, offering a further step toward understanding the gap between its worst-case bounds and empirical performance.

Probabilistic analyses of LP algorithms beyond simplex and interior-point methods are very limited. \citet{blum2002smoothed} show that the perceptron algorithm, a simple greedy approach, achieves high-probability polynomial smoothed complexity $\widetilde{O}(\frac{d^3m^2}{\delta^2})$ for finding feasible solutions with probability at least $1-\delta$ for LP instances with $m$ constraints and $d$ variables. 

There is some recent work on average-case analyses of first-order methods for unconstrained quadratic optimization problems, such as \citet{pedregosa2020acceleration,scieur2020universal,paquette2021sgd,paquette2023halting}. More recently, \citet{anagnostides2024convergence} study the high-probability performance of some first-order methods for zero-sum matrix games to eliminate condition number dependence.

\subsection{Notation}

For any positive integer $n$, let $[n]$ denote the set $\{1,2,\dots,n\}$. For a matrix $A \in \mathbb{R}^{m \times n}$, we denote by $A_{\cdot,i}$ its $i$-th column ($i \in [n]$) and by $A_{j,\cdot}$ its $j$-th row ($j \in [m]$). For any subset $\Theta$ of $[n]$, $A_{\Theta}$ denotes the submatrix formed by columns indexed by $\Theta$. 
For a vector $v$, $\|v\|$ denotes the Euclidean norm, and $\|v\|_1$ denotes the $\ell_1$ norm. For a matrix $A$, $\|A\|$ denotes the spectral norm and $\|A\|_F$ denotes the Frobenius norm. For scalars $a$ and $b$, we use $a \wedge b$ and $a \vee b$ to denote their minimum and maximum, respectively. We denote the singular values of $A\in \R^{m\times n}$ with $m\le n$ by $\sigma_1(A),\sigma_2(A),\dots,\sigma_{m}(A)$, where $\sigma_1(A)\ge \sigma_2(A)\ge \dots \ge \sigma_{m}(A)\ge 0$. The notation $0_m$ denotes the $m$-dimensional zero vector. 
Throughout this paper, $O(\cdot)$ denotes upper bounds, and $\Omega(\cdot)$ denotes lower bounds, both up to absolute constant factors if not specified. Similarly, $\widetilde{O}(\cdot)$ and $\widetilde{\Omega}(\cdot)$ denote upper and lower bounds respectively, while hiding additional logarithmic factors. For an event $E$, we use $E^c$ to denote the complement of $E$.

\section{Restarted PDHG for Linear Programming}\label{sec:preliminary}
Recall that this paper studies the standard-form LP problem \eqref{pro:primal LP}. When rows of $A$ are linearly independent, a basis of $A$ contains $m$ columns. We denote by $d := n - m$ the number of nonbasic columns, which is also the difference between $n$ and $m$. 
The corresponding dual problem of \eqref{pro:primal LP} is \eqref{pro:dual LP}. 
A more symmetric form of the dual problem can be obtained by eliminating $y$. Let $Q$ be any matrix so that the null space of $Q$ is equal to the image space of   $A^\top$, and let $\hat{x}$ be any feasible primal solution. Then  \eqref{pro:dual LP} can be equivalently reformulated in terms of $s$ alone as follows:
\begin{equation}\label{pro:dual LP_s}
	\min _{s \in \mathbb{R}^n} \ \hat{x}^\top s \quad \text { s.t. } Qs = Qc, s \geq 0 .
\end{equation}
This reformulation of the dual was first proposed by \citet{todd1990centered}. The optimal slack $s^\star$ of \eqref{pro:dual LP} is identical to the optimal solution $s^\star$ of \eqref{pro:dual LP_s}. With $s^\star$, any $(y,s^\star)$ such that $A^\top y + s^\star =c$ forms an optimal solution for \eqref{pro:dual LP}. Similarly, if $y^\star$ is an optimal dual solution, then $(y^\star, c - A^\top y^\star)$ is optimal for \eqref{pro:dual LP} and $c - A^\top y^\star$ is optimal for \eqref{pro:dual LP_s}.

The optimality conditions state that $x^\star$ and $s^\star$ are optimal if and only if they are feasible for \eqref{pro:primal LP} and \eqref{pro:dual LP_s} and satisfy the complementary slackness condition $(x^\star)^\top s^\star = 0$. 
Since $\hat{x}^\top s = \hat{x}^\top (c - A^\top y) = \hat{x}^\top c - b^\top y$ for any feasible solution $\hat{x}$ of \eqref{pro:primal LP} and feasible $(y,s)$ of \eqref{pro:dual LP}, the objective vector $\hat{x}$ of \eqref{pro:dual LP_s} can be replaced by any primal feasible solution $\hat{x}$ without altering the optimal solutions $x^\star$ and $s^\star$. Likewise, $x^\star$ and $s^\star$ remain optimal if the objective vector $c$ of \eqref{pro:primal LP} is replaced with any $\hat{c}$ satisfying $Q\hat{c} = Qc$ (namely, $\hat{c}=c + A^\top y_0$ for some $y_0 \in \R^m$). In later sections, we frequently assume $c$ has been replaced by $\bar{c} := c + A^\top \bar{y}$ during presolving, where $\bar{y} := \argmin_{y} \|c + A^\top y\|$. This substitution simplifies the analysis without affecting the optimal solutions $(x^\star, s^\star)$.  

Furthermore, the optimal solutions of primal problem \eqref{pro:primal LP} and dual problem \eqref{pro:dual LP} form the saddle point of the Lagrangian $L(x,y)$ that is defined as follows:
\begin{equation}\label{pro:saddlepoint_LP}
	\min_{x\in \R^n_+} \max_{y\in \R^m} \quad L(x,y) := c^\top x + b^\top y - (Ax)^\top y \ .
\end{equation} 
Any $x^\star$ and $(y^\star, c - A^\top y^\star)$ are optimal solutions of \eqref{pro:primal LP} and \eqref{pro:dual LP} if and only if $(x^\star, y^\star)$ is a saddle point of \eqref{pro:saddlepoint_LP}, and vice versa. The formulation \eqref{pro:saddlepoint_LP} is the problem that the primal-dual hybrid gradient method directly addresses. 

\subsection{Restarted Primal-dual hybrid gradient method (rPDHG)}\label{subsec:rPDHG}
A single iteration of PDHG, denoted by $(x^+, y^+) \leftarrow \textsc{OnePDHG}(x, y)$, is defined as follows:
\begin{equation}\label{eq_alg: one PDHG}
	\left\{\begin{array}{l}
		x^+ \leftarrow P_{\R^n_+}\left(x-\tau\left(c-A^{\top} y\right)\right) \\
		y^+ \leftarrow y+\sigma\left(b-A\left(2 x^+-x\right)\right)
	\end{array}\right. 
\end{equation}
where $P_{\R^n_+}\left(\cdot\right)$ denotes the projection onto $\R^n_+$, which means taking the positive part of a vector, and $\tau$ and $\sigma$ are the primal and dual step sizes, respectively. The \textsc{OnePDHG} iteration involves only two matrix-vector multiplications, one projection onto $\R^n_+$ and several vector-vector products. These operations avoid the computationally expensive matrix factorizations. Therefore, PDHG is well-suited for solving large LP instances and exploiting parallel implementation on modern computational architectures such as GPUs.

Furthermore, \citet{applegate2023faster} introduce the use of restarts to accelerate the convergence rate of PDHG. Here, the ``restarts'' mean that the method occasionally restarts from the average iterate of the previous consecutive many iterates. We will refer to this scheme as ``rPDHG'' for the restarted PDHG. Compared with the vanilla version, rPDHG can achieve linear convergence on LP instances and has shown strong practical performance. The rPDHG, together with some heuristic techniques, is the base method used in most current state-of-the-art first-order LP solvers.  

In this paper, we organize PDHG iterations into outer loops, and restart at the averaged iterate of the current inner loop.
Following \citet{applegate2023faster}, we use the $\beta$-restart criterion based on the \emph{normalized duality gap} (see, e.g., (4a) of \citet{applegate2023faster}): an outer loop restarts once the
normalized duality gap of the averaged iterate decreases to a $\beta$ fraction of its value at the start of that outer loop.
Throughout this paper (both theory and experiments), we take $\beta$ as a fixed absolute constant (we use $\beta=1/e$).
See Algorithm \ref{alg: PDHG with restarts} in Appendix \ref{appsec:restartPDHG} for the complete algorithm framework considered in this paper. 
Throughout the remainder of the paper, rPDHG refers to Algorithm \ref{alg: PDHG with restarts} in Appendix \ref{appsec:restartPDHG}.

Other restart triggers have also been used in practice. For example, \citet{applegate2023faster} discuss fixed-frequency restarts in addition to the normalized duality gap restart, and \citet{lu2024restarted} consider alternative triggers based on progress measures. These restart schemes have similar theoretical iteration bounds and empirical performance. No single trigger is consistently better than the others. We use the restart scheme based on the normalized duality gap in this paper, as it is widely used in implementations such as \citet{applegate2021practical,lu2023cupdlpnew}.

Note that rPDHG contains nested loops, but for notational simplicity, we use a single superscript indexed by the number of \textsc{OnePDHG} steps performed. We denote the primal and dual solutions of rPDHG after $k$ steps of \textsc{OnePDHG} by $x^k$ and $y^k$. The corresponding slack $c-A^\top y^k$ is denoted by $s^k$.

\section{High-Probability Computational Guarantees of rPDHG}\label{sec:high-probability} 
In this section, we describe the probabilistic model used in our analysis and present the polynomial-time complexity of rPDHG when applied to this model in a high-probability sense.

\subsection{Probabilistic linear programming model}\label{subsec:randomLPmodel}

Our probabilistic model builds on the classic model proposed by \citet{todd1991probabilistic}, in which the constraint matrix is drawn from a specific probability distribution, and the right-hand side and the objective vector are computed based on given random primal and dual solutions. Different versions of this model have been extensively studied in the average-case complexity analyses of interior-point methods (see, for example, \citet{ye1994toward,anstreicher1993averagenew,anstreicher1999probabilistic,huang2005expected,huang2000average,ji2008probabilistic}). 

Compared with other probabilistic models, a significant difference lies in the distribution of the constraint matrix. It is usually assumed that each component of the constraint matrix obeys a Gaussian distribution (see, e.g., \citet{anstreicher1999probabilistic} and \citet{todd1991probabilistic}). Our probabilistic model weakens this assumption to the more general case of sub-Gaussian distributions. Below, we provide the definition of sub-Gaussian random variables along with some typical examples (see proofs and additional examples in \citet{wainwright2019high}):
\begin{definition}[Sub-Gaussian random variable]
		A random variable $X$ is sub-Gaussian if there exists $\sigma >0$ such that $\E\left[e^{\lambda(X-\E[X])}\right]\le e^{\frac{\sigma^2\lambda^2}{2}}$ for all $\lambda \in \R$. The parameter $\sigma$ is referred to as the sub-Gaussian parameter.
\end{definition}
\begin{example}[Gaussian random variable]
	Let $X$ be a Gaussian random variable with mean $\mu$ and variance $\sigma^2$. Then $X$ is sub-Gaussian with the parameter $\sigma$.
\end{example}
\begin{example}[Bounded random variable]
		Let $X$ be a random variable supported on the bounded interval $[a,b]$. Then $X$ is sub-Gaussian with a parameter $\sigma$ that is less than or equal to $\frac{b-a}{2}$.
\end{example}

The family of sub-Gaussian distributions contains many commonly used distributions. In particular, it contains any bounded distribution, such as the Bernoulli distribution and the uniform distribution on a closed interval. 

Our probabilistic model involves the notions of sub-Gaussian matrices and nonnegative absolutely continuous sub-Gaussian vectors defined as follows.
\begin{definition}[Sub-Gaussian matrix]\label{def:random matrix}
	A matrix $A$ is a sub-Gaussian matrix if its elements are independent and identically distributed (i.i.d.), each obeying a mean-zero, unit-variance sub-Gaussian distribution.
	The sub-Gaussian parameter of each component is denoted by $\sigma_{A_{ij}}$, and the sub-Gaussian parameter of $A$ is defined to be $\sigma_A := \max_{ij}\sigma_{A_{ij}}$.
\end{definition} 
\begin{definition}[Nonnegative absolutely continuous sub-Gaussian vector]\label{def:positive_sub-Gaussian}
		A vector $u$ is a nonnegative absolutely continuous sub-Gaussian vector if its components are independent (and potentially different) nonnegative sub-Gaussian random variables whose probability density functions are bounded above by one. We denote the maximum of the means and the maximum of the sub-Gaussian parameters over all components of $u$ by $\mu_u$ and $\sigma_u$, respectively. 
\end{definition}

With the above definitions, we now describe the probabilistic model considered in this paper.
\begin{definition}[Probabilistic model]\label{def:random_lp}
	Instances of the probabilistic model are as follows:
	\begin{equation}\label{pro:RLP instance}
		\min_{x\in\R^n} \ \hat{s}^\top x \ , \quad \mathrm{s.t.} \ Ax = A\hat{x} \ , x \ge 0 \  .
	\end{equation}
	The constraint matrix $A$ is a sub-Gaussian matrix in $\R^{m \times n}$ as defined in Definition \ref{def:random matrix}, where $m < n$. The right-hand side and the objective vectors are computed based on $A$ and the primal and dual solutions $\hat{x}$ and $\hat{s}$, where $\hat{x}$ and $\hat{s}$ are generated as follows:
	\begin{equation}\label{eq:given_optimal_solution}
		\hat{x}=\begin{pmatrix}
			u^1 \\
			0_d
		\end{pmatrix} \text{ and }
		\hat{s}=\begin{pmatrix}
			0_m \\
			u^2
		\end{pmatrix}
	\end{equation}
	where $u^1 \in \R_+^m$ and $u^2 \in \R_+^{d}$, and the vector $u := (u^1, u^2)$ is a nonnegative absolutely continuous sub-Gaussian vector in $\R^n$ as defined in Definition \ref{def:positive_sub-Gaussian}. The random matrix $A$ and the random vector $u$ are independent.
 (The objective vector $\hat{s}$ can be replaced by any $c$ of the form $c=\hat{s}+A^\top y$ since the optimal solution sets $(\calX^\star,\calS^\star)$ remain unchanged.
	Optionally, $\hat{s}$ can be replaced by $\bar{c}$, the objective vector with the smallest norm, defined as $\bar{c} :=  \hat{s} + A^\top \hat{y}$ where $\hat{y} = \argmin_{y\in\R^m}\|\hat{s} + A^\top y\|$, and so satisfies $A\bar c = 0$.)
\end{definition} 

In the probabilistic model in Definition \ref{def:random_lp}, the components of $A$ have unit variance and the probability densities of the components of $u$ are at most one. If these requirements do not hold, they can be satisfied by scalar rescaling of the data instances.

	Instances of the probabilistic model have the option of using the objective vector $\bar{c}$ instead of $\hat{s}$.
As discussed above, this $\bar{c}$ satisfies $A\bar{c} = 0$ and does not change the optimal solution sets $(\calX^\star,\calS^\star)$ because $\bar{c} = \hat{s} + A^\top \hat{y}$.
Keeping $\hat{s}$ would introduce an additional minor term involving the component in the image space of $A^\top$ in the complexity bound (see Remark 3.5 in \citet{xiong2023computational}), which makes the bound look unnecessarily complicated.
Furthermore, computing $\bar{c}$ (e.g., using the conjugate gradient method to compute $\hat{y}$) is usually much less expensive than solving the LP itself, and $A\bar{c}=0$ is often assumed in rPDHG analyses such as \citet{xiong2024role,xiong2024accessiblenew}.
To enhance clarity and to align with these previously established results used in our analysis, we presume throughout this paper that the objective vector is set to $\bar{c}$ during a presolving step before applying rPDHG.

When applying rPDHG to instances of this model, we assume that rPDHG regards them as regular standard-form LP instances without knowing any prior information about the distribution of the input data.

The model in Definition \ref{def:random_lp} is analogous to the model of \citet{anstreicher1999probabilistic}, specifically the TDMV1 model. 
It is an important case of \citet{todd1991probabilistic}'s classic probabilistic model and has been studied by \citet{anstreicher1993averagenew,ye1994toward,ye1997interior} and others for analyzing the average-case performance of interior-point methods. Our model and those of \citet{anstreicher1999probabilistic,todd1991probabilistic} all assume that the constraint matrix and a pair of primal-dual feasible (and optimal) solutions are sampled from specific probability distributions, after which the right-hand side $b$ and the objective vector $c$ of the random LP instance are computed by $b = A\hat{x}$ and $c = \hat{s} $. 
One advantage of using this model is that any instance of this model is feasible and has a bounded optimal solution that is randomly distributed. 

\vspace{5pt}
\textbf{Distribution of the constraint matrix.} 
Our model allows the constraint matrix to follow a general sub-Gaussian distribution, whereas the probabilistic models most commonly studied in the LP literature assume Gaussian or other symmetry-based distributions. Generalizing the distributional assumptions is important because it broadens the range of random LP instances for which the analysis can offer insight into rPDHG's observed behavior. For example, in the probabilistic analysis of simplex methods, \citet{borgwardt1982average,smale1983average} assume that the polytopes come from a special spherically symmetric distribution. Later, \citet{adler1985simplex,adler1986family,todd1986polynomial} assume that the constraint matrix is drawn from a sign-invariant distribution.  \citet{spielman2004smoothed} study Gaussian perturbations of input data to the simplex method. On the other hand, most probabilistic analyses of interior-point methods have been built on \citet{todd1991probabilistic}'s probabilistic model, which assumes that the constraint matrix is a Gaussian matrix. Examples include \citet{mizuno1993adaptive,ye1994toward,anstreicher1993averagenew,anstreicher1999probabilistic}. 

We will show later in Section \ref{subsec:probabilistic_analysis_subgaussian} that, for our model with general sub-Gaussian distributions, rPDHG has polynomial-time complexity in a high-probability sense. Given the popularity of Gaussian distributions in the literature on probabilistic analysis, later in Section \ref{subsec:probabilistic_analysis_gaussian}, we will show that rPDHG has even better high-probability polynomial-time complexity when the constraint matrix is Gaussian.

\vspace{5pt}
\textbf{Distribution of solutions $\hat{x}$ and $\hat{s}$.}
In our probabilistic model, $\hat{x}$ and $\hat{s}$ are optimal primal-dual solutions because they are feasible primal and dual solutions and satisfy the complementary slackness condition.  Because rPDHG is invariant under permutations of variables, without loss of generality, in the model we assign the indices of the possible nonzero components of $\hat{x}$ and $\hat{s}$ to $\{1,\dots,m\}$ and $\{m+1,\dots,n\}$. 
We use $B$ and $N$ to denote the submatrices of $A$ corresponding to the column indices $\{1,\dots,m\}$ and $\{m+1,\dots,n\}$, respectively, so that $A$ is represented as $(B, N)$. When $B$ is full-rank and $(\hat{x},\hat{s})$ satisfies the strict complementary slackness condition, the optimal solution sets $\calX^\star$, $\calY^\star$ and $\calS^\star$ are all singletons. It is formally stated as follows (proof in Appendix \ref{appsec:high-probability}):
\begin{lemma}\label{lm:full-rank-unique-optima}
	For an instance of the probabilistic model, $\calX^\star$, $\calY^\star$ and $\calS^\star$ are all singletons with $\calX^\star=\{\hat{x}\}$ and $\calS^\star = \{\hat{s}\}$ if and only if $B$ is full-rank and $(\hat{x},\hat{s})$ satisfies the strict complementary slackness condition, namely, $u = (u^1,u^2) > 0$.
\end{lemma} 
\noindent
Furthermore, once $m$ is sufficiently large, it is highly probable that the instance of the probabilistic model has unique primal and dual optimal solutions.
\begin{lemma}\label{lm:full-rank-probability}
	There exists a constant $\nu>0$ that depends only (and at most polynomially) on the sub-Gaussian parameter of $A$, such that the probability of $B$ being full-rank is at least $1 - e^{-\nu m}$.
\end{lemma}
\noindent The proof of Lemma \ref{lm:full-rank-probability} is provided in Section \ref{subsec:helper_lemmas}.

The distributions of $\hat{x}$ and $\hat{s}$ also differ between our probabilistic model and that of \citet{anstreicher1999probabilistic}. On the one hand, the model of \citet{anstreicher1999probabilistic} permits varying numbers of nonzeros in $\hat{x}$ and $\hat{s}$ while maintaining strict complementarity. This may result in multiple optimal solutions with high probability. On the other hand, the components of $u$ in their model are i.i.d. from a folded Gaussian distribution (the distribution of the absolute value of a Gaussian random variable). In contrast, our model fixes the numbers of nonzeros in $\hat{x}$ and $\hat{s}$ to $m$ and $d$, respectively, while allowing the components of $u$ to follow potentially different sub-Gaussian distributions.

We acknowledge that for a very large number of LP instances occurring in practice, the optimal solution is not unique. But the unique optimum property is generic in theory, especially in probabilistic models, because most randomly generated LP problems (unless specially designed) are nondegenerate (see, for example,  \citet{ye1994toward,spielman2003smoothed,borgwardt1982some,todd1986polynomial,adler1986family}). 

The flexibility in the distribution of $u$ in our model provides tools for various aims of analyses. When $u$ is a random vector of the folded Gaussian distribution, the probabilistic model is for the average-case analysis. 
When $u$ is a given fixed vector with random perturbations, the probabilistic model is for smoothed analysis on the dependence of the optimal solution. Furthermore, we will show later in Section \ref{sec:effect_of_optimal_solution} that the performance of rPDHG is highly influenced by the distribution of the optimal solution, and the unique optimum property of our model is helpful in generating artificial LP instances of various difficulty levels.

It should be mentioned that \citet{todd1991probabilistic} also provides a model that allows control over the degree of degeneracy in $\hat{x}$ and $\hat{s}$ (TDMV2 of \citet{anstreicher1999probabilistic}). But Todd and \citet{anstreicher1999probabilistic} later pointed out a subtle error in the analysis of this model. Due to this error, several papers, including \citet{anstreicher1993averagenew} and \citet{ye1997interior}, that had claimed to analyze the average-case performance of interior-point methods on this model actually applied only to the TDMV1 model of \citet{anstreicher1999probabilistic}.  This subtle error also affected several papers that studied the version of \citet{todd1991probabilistic}'s probabilistic model with $\hat{x} =e$ and $\hat{s} = e$. These errors were later corrected by \citet{huang2000average} and \citet{ji2008probabilistic}.  

Next, in Section \ref{subsec:probabilistic_analysis_subgaussian}, we present our main results on the high-probability polynomial-time complexity of rPDHG for instances of the probabilistic model with sub-Gaussian input data. In Section \ref{subsec:probabilistic_analysis_gaussian}, we then consider instances with Gaussian input data.

\subsection{High-probability performance guarantees for LP instances with  sub-Gaussian input data} \label{subsec:probabilistic_analysis_subgaussian}

In this subsection, we analyze the performance of rPDHG on instances of the probabilistic model defined in Definition \ref{def:random_lp}. Our focus is on the dependence of the iteration bounds on the dimensions $m$ and $n$.  The bounds may also contain some constants that depend only on the parameters of the distributions of $A$ and $u$, namely $\sigma_A$, $\mu_u$ and $\sigma_u$ in Definitions \ref{def:random matrix} and \ref{def:positive_sub-Gaussian}.  
Below we define an $\eps$-optimal solution for the LP primal and dual problems \eqref{pro:primal LP} and \eqref{pro:dual LP}.
\begin{definition}[$\eps$-optimal solutions]
	A primal-dual pair $(x,s)$ is an $\eps$-optimal solution if there exists a primal-dual optimal pair $(x^\star,s^\star)$ such that $\|(x,s) - (x^\star,s^\star)\|\le \eps$.
\end{definition}

It is often observed in practice (and shown in theory) that rPDHG exhibits a ``two-stage phenomenon,'' in which its iterations can be divided into Stage I and Stage II. In Stage I, the iterates eventually reach the point where their positive components correspond to the optimal basis. In Stage II, rPDHG exhibits fast local convergence to an optimal solution (see, e.g., \citet{xiong2024accessiblenew,lu2024geometrynew,lu2024restarted}). In other words, there exists $T_{basis}$ such that, for all $t\ge T_{basis}$, $x^{t}_{i} > 0$ if and only if $i$ is in the optimal basis, after which rPDHG exhibits faster local linear convergence to an optimum. We let $T_{local}$ denote the number of additional iterations beyond $T_{basis}$ required to compute an $\eps$-optimal solution $(x^t,s^t)$. The following theorem presents our high-probability bounds on $T_{basis}$ and $T_{local}$.

\begin{theorem}\label{thm:complexity-random}
	Suppose that rPDHG is applied to an instance of the probabilistic model in Definition \ref{def:random_lp} (with objective vector $\bar{c}$). Let $c_0:=\mu_u+2\sigma_u$. There exist constants $C_0,C_1,C_2 > 0$ that depend only (and at most polynomially) on $\sigma_A$, for which the following high-probability iteration bounds hold:
	\begin{enumerate}
		\item (Optimal basis identification) 
		\begin{equation}\label{eq:optimal_basis_random}
			      \Pr\left[
			      T_{basis}
			      \le
			      \frac{m^{0.5}n^{3.5}}{d+1}\cdot\frac{  c_0C_1 \cdot \ln(3/\delta)}{\delta}  \cdot \ln\left( \tfrac{m^{0.5}n^{3.5}}{d+1}\cdot\tfrac{  c_0C_1 \cdot \ln(3/\delta)}{\delta} 
			      \right)
			      \right]\ge 1 - \delta - 6\left(\tfrac{1}{e^{C_0}}\right)^m - \left(\tfrac{1}{2}\right)^{d+1}
		      \end{equation}
		      for any $\delta \in (0,1]$.
		\item\text{(Fast local convergence)} 	
		Let $\eps > 0$ be any given tolerance.  
		\begin{equation}\label{eq:local_convergence_random}
			      \Pr\left[
			      T_{local}
			      \le
				  m^{0.5} n^{0.5} \cdot 
			      \frac{C_2  }{\delta}\cdot  \max\left\{0,\ln\left(\frac{c_0}{\eps}\right)\right\}
				      \right]\ge 1 - \delta - 3\left(\tfrac{1}{e^{C_0}}\right)^m
		      \end{equation}
		      for any $\delta \in (0,1]$.
	\end{enumerate} 
\end{theorem}

In the above theorem, we have divided the rPDHG iterations into two stages and presented probabilistic bounds for them separately. Technically speaking, once the iterates settle on the optimal basis, the optimal solution $(x^\star,s^\star)$ could be directly computed by solving two linear systems---one for the primal basic system and the other for the dual basic system. (This is common in finite termination methods for interior-point methods, such as \citet{ye1992finite}, where the critical effort lies in computing projections.) Indeed, one could compute optimal solutions in finite time using rPDHG by first running rPDHG for $T_{basis}$ \textsc{OnePDHG} iterations and then solving the two associated linear systems. However, this is not a practical approach for at least two reasons. First, determining whether $(x^t,s^t)$ has already settled on the optimal basis can be difficult. Second, rPDHG automatically exhibits fast local convergence after identifying the optimal basis, so computing an $\eps$-optimal solution requires at most $T_{basis}+T_{local}$ iterations. The Stage-II bound in \eqref{eq:local_convergence_random} depends only logarithmically on the target accuracy and on $c_0$. This result does not conflict with the global linear convergence guarantee established by \citet{applegate2023faster}. Actually, rPDHG converges linearly in Stage I. We discuss this further in Remark \ref{rmk:two-stage-linear-convergence} in Section \ref{sec:prove_thm_random}.

The inequalities \eqref{eq:optimal_basis_random} and \eqref{eq:local_convergence_random} are high-probability upper bounds for $T_{basis}$ and $T_{local}$. The right-hand sides of \eqref{eq:optimal_basis_random} and \eqref{eq:local_convergence_random} both contain $1-\delta$ as well as some additional terms that decrease exponentially in the dimensions $m$ and $d$. In other words, if $m$ and $d$ are sufficiently large -- essentially larger than $O(\ln\frac{1}{\delta})$ -- these additional terms are negligible, and so \eqref{eq:optimal_basis_random} and \eqref{eq:local_convergence_random} describe the upper bounds on $T_{basis}$ and $T_{local}$ with a probability that is roughly equal to $1-\delta$. 
% Because the goal of probabilistic analysis is mainly on revealing the dependence on dimensions when the dimensions are large,  
Similar requirements on the dimensions being sufficiently large are common in other probabilistic analyses of linear programming, such as \citet{mizuno1993adaptive,ye1994toward,huang2005expected} for interior-point methods, and \citet{borgwardt1977b,adler1986family,shamir1987efficiency} for simplex methods. 

The constant $c_0$ depends only on $\mu_u$ and $\sigma_u$, whereas $C_0,C_1$, and $C_2$ depend only on $\sigma_A$. See Remark \ref{rmk:constants-subgaussian} in Section \ref{sec:prove_thm_random} for further discussion of these constants. The term $\frac{m^{0.5}n^{3.5}}{d+1}$ is at most as large as $\frac{m^{0.5}n^{3.5}}{2}$.  When $d$ (recall $d:=n-m$) is not too small compared to $n$ in the sense that $d+1 \ge n/C_3$ for some absolute constant $C_3 >1$, then the term $\frac{m^{0.5}n^{3.5}}{d+1}$ is $O(m^{0.5}n^{2.5})$. 

The corollary below summarizes the high-probability iteration bounds when (i) $m$ is not too close to $n$, and (ii) $\delta$ is not exponentially small in $d$ and $m$: 

\begin{corollary}\label{cor:complexity-random}
	Let $C_3 > 1$ be a given absolute constant. Suppose that rPDHG is applied to an instance of the probabilistic model (with objective vector $\bar{c}$). When $d$ satisfies $\frac{n}{d+1}\le C_3$, it holds with probability at least $1-\delta$ that rPDHG computes an $\eps$-optimal solution within at most
	$$
		\widetilde{O}\left(
		\frac{ n^{2.5}m^{0.5} }{\delta} \right) + O\left( \frac{n^{0.5}m^{0.5}}{\delta}
		\cdot  \ln\left(\frac{1}{\eps}\right) 
		\right)
	$$
	iterations for all $\eps \in (0,e^{-1}]$ and any $\delta \in (0,1]$ satisfying $\delta > 11 \cdot \max\Big\{\left(\tfrac{1}{e^{C_0}}\right)^m , \left(\tfrac{1}{2}\right)^{d+1}\Big\}$.
	Moreover, $T_{basis}$ is at most $\widetilde{O}\Big(\frac{n^{2.5}m^{0.5}}{\delta} \Big)$.   Here $\widetilde{O}(\cdot)$ omits factors of an absolute constant, $c_0,C_1,C_3$ and logarithmic terms that 
	involve $c_0,C_1,C_3,m,n$ and $\frac{1}{\delta}$, and $O(\cdot)$ additionally omits factors of  $C_2$ and $\ln c_0$.
\end{corollary}

Corollary \ref{cor:complexity-random} shows that under a mild assumption on $\delta$ and on the dimensions $m$ and $d$ (the interval for $\delta$ is nonempty only when $m > \tfrac{\ln(11)}{C_0}$ and $d > \log_2(11)-1\approx2.4594$), rPDHG settles on the optimal basis within  $\widetilde{O}\Big(\frac{n^{2.5}m^{0.5}}{\delta} \Big)$ iterations, and it computes an $\eps$-optimal solution in an additional $O\Big(\frac{n^{0.5}m^{0.5}}{\delta} \cdot  \ln\frac{1}{\eps} \Big)$ iterations, with probability at least $1-\delta$. 
To the best of our knowledge, the above result is the first high-probability polynomial-time complexity bound for rPDHG under a probabilistic input model. It provides a probabilistic explanation for why rPDHG can avoid its exponentially poor worst-case complexity on a broad class of random LP instances, thereby making progress toward understanding the gap between its worst-case theory and observed performance. Moreover, it shows that the dependence of the high-probability complexity on $n$ in Stage I (settling on the optimal basis) is higher than in Stage II (local convergence). This observation also aligns with the two-stage phenomenon reported in worst-case complexity and practical experimental results (see \citet{xiong2024accessiblenew} and \citet{lu2024restarted,lu2024geometrynew}).

Furthermore, to the best of our knowledge, this is the first probabilistic iteration complexity bound for an LP algorithm that allows the constraint matrix to have general, not necessarily Gaussian, sub-Gaussian entries. The Gaussian matrix is relatively easy to analyze because it has nice symmetry in the sense of geometry, its range is the unique orthogonally invariant distribution, and the corresponding LP instance is nondegenerate with probability 1 (see \citet{todd1991probabilistic}). The more general sub-Gaussian random matrix model is harder to analyze than it may look. Even the behavior of the extreme singular values of sub-Gaussian matrices was not well understood until about 15 years ago (see the invited lecture at ICM 2010 by \citet{rudelson2010non}). Our analysis approach relies on the nonasymptotic result of \citet{rudelson2009smallest} on the smallest singular value of a sub-Gaussian matrix.  Proofs of our results above are presented in Section \ref{sec:prove_thm_random}.

It should be noted that real-world LP instances often differ from our model by having multiple optimal solutions or input data that are not well-approximated by sub-Gaussian distributions (e.g., heavy-tailed or highly structured data). In such cases, we do not yet have high-probability performance guarantees. 
The two-stage characterization used in Theorem \ref{thm:complexity-random} is not well-defined when the optimal basis is not unique, and our probabilistic analysis relies on the accessible iteration bound of \citet{xiong2024accessiblenew}, which is for LP instances with unique optimal solutions.  A possible extension to general degenerate instances is to start from the global convergence analyses that apply to general LPs, such as \citep{applegate2023faster} and \citep{xiong2023computational}, and then study the typical size of the condition measures used in these results under the probabilistic model.  Moreover, sub-Gaussianity enables us to invoke nonasymptotic random matrix singular value estimates. Extending the analysis to heavy-tailed or structured inputs would require different tools (e.g., \citet{cook2018lower,dumitriu2024extreme}). These are left for future investigation.

\subsection{High-probability performance guarantees for LP instances with Gaussian input data}\label{subsec:probabilistic_analysis_gaussian}

In this subsection, we analyze the performance of rPDHG under the probabilistic model wherein the constraint matrix $A$ is a Gaussian matrix, which we now define.
\begin{definition}[Gaussian matrix]\label{def:gaussian_matrix}
	A matrix $A$ is called a Gaussian matrix if its elements are i.i.d., each obeying the mean-zero Gaussian distribution with unit variance. 
\end{definition} 

A Gaussian matrix $A$ is a special case of a sub-Gaussian matrix in Definition \ref{def:random matrix}, with its sub-Gaussian parameter $\sigma_A$ equal to $1$ (see \citet{wainwright2019high}). In this section, we show that the high-probability iteration bound of rPDHG, particularly for Stage I (identifying the optimal basis), can be further improved in this classical Gaussian setting.  Similar to Theorem \ref{thm:complexity-random}, the theorem below presents high-probability bounds on $T_{basis}$ and $T_{local}$ when the matrix $A$ is Gaussian.

\begin{theorem}\label{thm:complexity-gaussian}
	Suppose that rPDHG is applied to an instance of the probabilistic model in Definition \ref{def:random_lp} (with objective vector $\bar{c}$) and suppose that the constraint matrix is a Gaussian matrix. Let $c_0 := \mu_u + 2\sigma_u$. There exist absolute constants $C_0,C_1,C_2>0$ for which the following high-probability iteration bounds hold:
	\begin{enumerate}
		\item\text{(Optimal basis identification)} 
		\begin{equation}\label{eq:optimal_basis_gaussian} 
			\begin{aligned}
			& \Pr\left[
			T_{basis}
			\le \frac{n^{3.5}}{d+1} \cdot \frac{c_0 C_1 \cdot \ln(6/\delta)\sqrt{\ln(4n/\delta)}}{\delta}
			\cdot \ln\left(
    \tfrac{n^{3.5}}{d+1} \cdot \tfrac{c_0 C_1 \cdot \ln(6/\delta)\sqrt{\ln(4n/\delta)}}{\delta}
\right)			\right]
			\\
			 &   \ge 1 - \delta -   (n+5)\left(\tfrac{1}{e^{C_0}}\right)^{m \wedge d}
			\end{aligned}
		\end{equation}
		      for any $\delta \in (0,1]$.
		\item\text{(Fast local convergence)} 	Let $\eps > 0$ be any given tolerance. 
		\begin{equation}\label{eq:local_convergence_gaussian}
			\Pr\left[
			T_{local}
			\le
			m^{0.5} n^{0.5}\cdot  \frac{C_2  }{\delta}
			\cdot \max\left\{0,\ln\left(\frac{c_0}{\eps}\right)\right\}
			\right]\ge 1 - \delta - 3\left(\tfrac{1}{e^{C_0}}\right)^m
		\end{equation}
		for any $\delta \in (0,1]$.
	\end{enumerate}
\end{theorem}

The inequalities \eqref{eq:optimal_basis_gaussian} and \eqref{eq:local_convergence_gaussian} are high-probability upper bounds on $T_{basis}$ and $T_{local}$. When $d$ and $m$ are sufficiently large relative to $\frac{n}{\delta}$, specifically when $m \wedge d \ge \Omega(\ln\frac{n}{\delta})$, the right-hand side of each inequality approaches $1-\delta$.  However, in the extreme case wherein $m$ or $d$ is too small compared with $n$, these bounds become trivial for all $\delta \in (0,1]$ since their right-hand sides become nonpositive.

Notice in Theorem \ref{thm:complexity-gaussian} that $C_0,C_1$, and $C_2$ are absolute constants. Indeed, these constants otherwise depend only on $\sigma_A$, and $\sigma_A=1$ for a Gaussian matrix. See Remark \ref{rmk:constants-gaussian} in Section \ref{sec:refined_analysis_gaussian} for further discussion of these constants.
The term $\frac{n^{3.5}}{d+1}$ is at most as large as $\frac{n^{3.5}}{2}$.  When $d$ is not too small compared to $n$, $\frac{n^{3.5}}{d+1}$ is $O(n^{2.5})$. Similarly to Corollary \ref{cor:complexity-random}, the corollary below summarizes the high-probability iteration bounds when (i) $m$ is not too close to $n$, and (ii) $\delta$ is not exponentially small:

\begin{corollary}\label{cor:complexity-gaussian}
	Let $C_3 > 0$ be any given absolute constant. Suppose that rPDHG is applied to an instance of the probabilistic model (with objective vector $\bar{c}$), and suppose that the constraint matrix $A$ is a Gaussian matrix. If $d$ satisfies $\frac{n}{d+1} \le C_3$, then it holds with probability at least $1 - \delta$ that rPDHG computes an $\eps$-optimal solution within at most
	$$
		\widetilde{O}\left(
		\frac{ n^{2.5}  }{\delta} \right) + 
		O\left(\frac{n^{0.5}m^{0.5}}{\delta}
		\cdot  \ln\left(\frac{1}{\eps}\right) 
		\right)
	$$ 
	iterations, for all $\eps \in (0,e^{-1}]$ and any $\delta \in (0,1]$ satisfying $\delta > 2(n+5)e^{-C_0(m\wedge d)}$.	
	Moreover, $T_{basis}$ is at most $\widetilde{O}\Big(\frac{n^{2.5}}{\delta} \Big)$. Here $\widetilde{O}(\cdot)$ omits factors of constants $c_0:= \mu_u+2\sigma_u$, an absolute constant and logarithmic terms that involve an absolute constant, $m, n$ and $ \frac{1}{\delta}$, and the $O(\cdot)$ additionally omits factors of  $C_2$ and $ \ln c_0$.
\end{corollary}

Corollary \ref{cor:complexity-gaussian} shows that under a mild assumption on $\delta$ and the dimensions $m$ and $d = n-m$ (the interval for $\delta$ is nonempty only when $m\wedge d > \frac{1}{C_0}\ln\!\big(2(n+5)\big)$), rPDHG settles on the optimal basis within  $\widetilde{O}\Big(\frac{n^{2.5}}{\delta} \Big)$ iterations, and it computes an $\eps$-optimal solution in an additional $O\Big(\frac{n^{0.5}m^{0.5}}{\delta} \cdot  \ln\frac{1}{\eps} \Big)$ iterations, with probability at least $1-\delta$. 
Comparing Corollary \ref{cor:complexity-gaussian} with Corollary \ref{cor:complexity-random} (for sub-Gaussian matrices), the dependence on $m^{0.5}$ is eliminated in the Stage-I iteration bound while the Stage-II iteration bound remains almost identical.

Now we compare this high-probability iteration bound with the probabilistic bounds of interior-point methods and simplex methods. In general, the interior-point method has far less dependence on the dimension in the number of iterations, but the per-iteration complexity is significantly higher than that of rPDHG because it needs to solve a normal equation of a dense normal matrix. The model of \citet{anstreicher1999probabilistic} is the most similar one to ours, and \citet{anstreicher1999probabilistic} prove that the expected number of iterations of an interior-point method is at most $O(n\cdot \ln(n))$. They also remark that their proof can be easily modified for other interior-point methods and derive an $O(n^2\ln(n))$ expected iteration bound for the interior-point methods of \citet{zhang1994convergence,wright1994infeasible,zhang1994superlinear} and others. 
\citet{huang2005expected} shows that the expected and high-probability numbers of iterations are both bounded above by $O(n^{1.5})$ in another model of \citet{todd1991probabilistic}.
\citet{ye1994toward} proves that the number of iterations is at most $O(\sqrt{n}\cdot\ln\frac{1}{\eps})$ with high probability under a different model.

As for simplex methods, direct comparison of the probabilistic bounds between rPDHG and simplex methods is challenging because different models and forms of LP are used in probabilistic analyses. Nevertheless, our high-probability iteration bound for rPDHG demonstrates comparable polynomial-time complexity to those established for simplex methods, with the additional advantage that rPDHG requires only two matrix-vector multiplications per iteration. To make the comparison as fair as possible, we consider problems where the number of constraints is of the same order as the number of variables, which we refer to as the dimension of the problem. \citet{todd1986polynomial,adler1985simplex,adler1987simplex} prove that the expected number of steps of several equivalent simplex methods on a certain probabilistic model is bounded by a quadratic function of the dimension.  \citet{adler1985simplex} also prove that the dependence is tight.  \citet{spielman2004smoothed} use the smoothed analysis framework to study models with Gaussian perturbations and derive a polynomial bound on the expected number of simplex pivots. \citet{dadush2018friendly} significantly simplify the analysis and establish a tighter polynomial bound with dimension exponent at most $3.5$. 

We note that Theorem \ref{thm:complexity-gaussian} does not imply a bound on the expected number of iterations of rPDHG. Indeed, the tails of \eqref{eq:optimal_basis_gaussian} and \eqref{eq:local_convergence_gaussian} do not decay to zero. 
In probabilistic analyses of linear programming, high-probability complexity bounds are also established by \citet{ye1994toward,blum2002smoothed} and others.  We are unable to prove bounds on the expected number of iterations for rPDHG, which is in contrast to the case of interior-point methods and simplex methods where such results are proven in \citet{anstreicher1999probabilistic,todd1986polynomial,dadush2018friendly} and others.

\section{Proof of Theorem \ref{thm:complexity-random}}\label{sec:prove_thm_random}

In this section, we prove Theorem \ref{thm:complexity-random}.
Section \ref{subsec:helper_lemmas} introduces several useful helper lemmas frequently used in the proofs, including concentration inequalities for sub-Gaussian random variables, bounds on the extreme singular values of random matrices, and a tail bound for the product of two independent random variables. Section \ref{subsec:worst-case} recalls the worst-case iteration bound of rPDHG. Section \ref{subsec:proof_of_probabilistic_analysis} contains the detailed proof of Theorem \ref{thm:complexity-random}.
 
\subsection{Lemmas of random variables and random matrices}\label{subsec:helper_lemmas}
The Hoeffding bound provides a tail bound for the sum of independent sub-Gaussian random variables (see \citet{wainwright2019high}).

\begin{lemma}[Hoeffding bound]\label{lm:hoeffding}
	Suppose that the sub-Gaussian variables $\{X_i\}_{i=1}^n$ are independent, and each $X_i$ has the sub-Gaussian parameter $\sigma_i$. Then for all $t \geq 0$, we have
	$\Pr\left[\sum_{i=1}^n\left(X_i-\E[X_i]\right) \geq t\right] \leq \exp \left(-\tfrac{t^2}{2 \sum_{i=1}^n \sigma_i^2}\right)$.
\end{lemma}

For a random matrix as defined in Definition \ref{def:random matrix}, the following result from \citet{rudelson2009smallest} provides a tail bound for the smallest singular value of the random matrix.
\begin{lemma}[Theorem 1.1 of \citet{rudelson2009smallest}]\label{lm:min_singular_value_rectangular}
	Let $A$ be a random matrix as defined in Definition \ref{def:random matrix}, where $A \in \R^{m \times n}$ and $n \ge m$. For every $\varepsilon > 0$, we have
	\begin{equation}\label{eq:lm:min_singular_value_rectangular}
		\Pr\left[\sigma_{m}(A) \leq \varepsilon\left(\sqrt{n}-\sqrt{m-1}\right)\right] \leq(C_{rv1} \varepsilon)^{d+1}+e^{-C_{rv2} n}
	\end{equation}
	where $C_{rv1}, C_{rv2}>0$ depend only (and at most polynomially) on the sub-Gaussian parameter $\sigma_A$. 
\end{lemma}
\noindent 
It is worth noting that the smallest singular value may not always be strictly greater than zero, which implies that the matrix $A$ may not be full-rank. However, the above lemma shows that as $n$ increases, the probability of $A$ being full-rank becomes very high. A direct application of this result is the proof of Lemma \ref{lm:full-rank-probability}.

\begin{proof}{Proof of Lemma \ref{lm:full-rank-probability}.}
	According to Lemma \ref{lm:full-rank-unique-optima}, the probability of Condition \ref{cond:unique_optima} being true is equal to $\Pr\left[\sigma_{m}(B) > 0\right]$. Since $B$ is a random matrix as defined in Definition \ref{def:random matrix} with $B \in \R^{m \times m}$, Lemma \ref{lm:min_singular_value_rectangular} implies that $\operatorname{Pr}\left[\sigma_m(B)=0\right]=\operatorname{Pr}\left[\sigma_m(B) \leq 0\right] \leq \lim _{\varepsilon \downarrow 0} \operatorname{Pr}\left[\sigma_m(B) \leq \varepsilon(\sqrt{m}-\sqrt{m-1})\right] \leq \lim _{\varepsilon \downarrow 0}\left(\left(C_{r v 1} \varepsilon\right)^1+e^{-C_{r v 2} m}\right)=e^{-C_{r v 2} m}.$	This completes the proof.
\hfill\Halmos\end{proof}

The largest singular value of random matrices is known to be upper bounded by $O(\sqrt{n})$ with high probability. It is formally stated in the following lemma.
\begin{lemma}\label{lm:max_singular_value_rectangular}
	Let $A$ be a random matrix as defined in Definition \ref{def:random matrix}, where $A \in \R^{m \times n}$ and $n \ge m$. Then,
	$$
		\Pr\left[\sigma_{1}(A) \geq 10\sigma_A\sqrt{n} \right]  \le e^{-6n} \ .
	$$ 
\end{lemma}

Next, we prove a tail bound of the product of two independent random variables with known heavy tails. 
\begin{lemma}\label{lm:tail_of_product}
	Let $X$ and $Y$ be two independent nonnegative random variables, and suppose that there exist $C_1, C_2, \delta_1,\delta_2 > 0$ such that for all $T > 0$:
	\begin{equation}\label{eq:tail_of_product_1}
		\Pr[X \ge T] \le \frac{C_1}{T} + \delta_1 \quad \text{ and }\quad
		\Pr[Y \ge T] \le \frac{C_2}{T} + \delta_2 \ .
	\end{equation}
	Then for any $\delta \in (0,1]$, the following inequality holds:
	\begin{equation}\label{eq:tail_of_product_1_2}
		\Pr\left[XY \ge \frac{6C_1C_2 \cdot \ln(3/\delta)}{\delta}\right] \le \delta + \delta_1 + 2 \delta_2 \ .
	\end{equation}
\end{lemma}

Proofs of Lemmas \ref{lm:max_singular_value_rectangular} and \ref{lm:tail_of_product} are provided in Appendix \ref{appsec:prove_thm_random}.

\subsection{Worst-case iteration bounds of rPDHG}\label{subsec:worst-case}

According to Lemmas \ref{lm:full-rank-unique-optima} and \ref{lm:full-rank-probability}, for instances of the probabilistic model, the probability of having full-rank $B$ and unique optimal solution $\calX^\star = \{\hat{x}\}$ and $\calS^\star = \{\hat{s}\}$ is at least $1 - e^{-\nu m}$. It happens with high probability when $m$ is sufficiently large. Therefore, in this subsection, we recall a theoretical iteration bound under the following condition of a unique optimum.
\begin{condition}\label{cond:unique_optima}
	The problem \eqref{pro:primal LP} has a unique optimal solution $x^\star$, and the dual problem \eqref{pro:dual LP} has a unique optimal solution $(y^\star,s^\star)$.
\end{condition}

Recently, \citet{xiong2024accessiblenew} proves an accessible iteration bound for rPDHG applied to LPs under Condition \ref{cond:unique_optima}. This new iteration bound is in closed form of the optimal solution and optimal basis.
Suppose that the optimal basis of $x^\star$ is $\{1,2,\dots,m\}$. We still let $B$ and $N$ denote the submatrices consisting of the columns indexed by $\{1,2,\dots,m\}$ and $\{m+1,m+2,\dots,n\}$. By Lemma \ref{lm:full-rank-unique-optima}, $B$ must be invertible. Then the iteration bound relies on the key quantity $\simplephi$ defined as follows:
\begin{equation}\label{def_geometric_measure}
	\simplephi :=  \big(\|x^\star+s^\star\|_1\big) \cdot \max\left\{
	\max_{1 \le j \le d} \frac{\sqrt{\left\|(B^{-1} N)_{\cdot,j}\right\|^2+1} }{s^\star_{m+j}} \, , \ \max_{1\le i \le m} \frac{\sqrt{\left\|(B^{-1} N)_{i,\cdot}\right\|^2+1} }{x^\star_i}  
	\right\}\ . 
\end{equation}
If Condition \ref{cond:unique_optima} does not hold, for notational simplicity, we let $\simplephi := \infty$. Then we have the following worst-case iteration bounds for $T_{basis}$ and $T_{local}$:

\begin{theorem}[Theorem 4.1 of \citet{xiong2024accessiblenew}]\label{thm:closed-form-complexity}
	Suppose that $Ac = 0$ and rPDHG is applied to solve the LP instance \eqref{pro:primal LP}. The following iteration bounds hold:
	\begin{enumerate} 
		\item\text{(Optimal basis identification)} There exists an absolute constant $\check{c}_1 > 1$ such that:
		      \begin{equation}\label{eq:basis_identification_T}
			       T_{basis} \le \check{c}_1 \cdot \kappa  {\simplephi} \cdot   \ln\left(\kappa  {\simplephi}  \right)  \ .
		      \end{equation}
		\item\text{(Fast local convergence)} There exists an absolute constant $\check{c}_2 > 0 $ such that:
		      \begin{equation}\label{eq:local_linear_convergence_T}
			    T_{local} \le \check{c}_2 \cdot \|B^{-1}\|\|A\|    \cdot \max\left\{0,\ \ln\left( \frac{\min_{1\le i \le n}\left\{x_i^\star + s_i^\star\right\}}{\eps}\right)\right\} .
		      \end{equation}
	\end{enumerate} 
\end{theorem}

Besides the ``complex'' expression of $\simplephi$ in \eqref{def_geometric_measure}, $\simplephi$ has the following simpler upper bound:
\begin{lemma}[Proposition 3.1 of \citet{xiong2024accessiblenew}]\label{lm:upperbound_phi}
	When Condition \ref{cond:unique_optima} holds, 
	$
		\simplephi    \le \frac{ \|x^\star+s^\star\|_1}{\min_{1\le i \le n}\left\{x_i^\star + s_i^\star\right\}} \cdot \|B^{-1} A \| .
	$
\end{lemma}

	Furthermore, it is actually indicated by \citet{xiong2024accessiblenew} that $\Phi$ is equal to $\frac{\|x^\star + s^\star\|_1}{\zeta_p \wedge \zeta_d}$ in which $\zeta_p$ and $\zeta_d$ are equivalent to three types of condition measures for primal and dual problems respectively. They are (i) stability under data perturbations, (ii) proximity to multiple optima, and (iii) the LP sharpness of the instance.

\subsection{Proof of Theorem \ref{thm:complexity-random}}\label{subsec:proof_of_probabilistic_analysis}

In this section, we prove Theorem \ref{thm:complexity-random}. Let $B$ be the $m \times m$ submatrix formed by the first $m$ columns of $A$. We denote the condition number of $A$ by $\kappa$, defined as $\kappa := \frac{\sigma_1(A)}{\sigma_m(A)}$. In addition, we let $\varphi$ denote the quantity defined as follows:
\begin{equation}\label{eq:def_varphi}
\varphi:=\frac{\|\hat{x} + \hat{s}\|_1}{\min_{1 \le i \le n} (\hat{x}_i + \hat{s}_i)} \ .
\end{equation}
We set $\varphi := \infty$ when $\min_{1 \le i \le n} (\hat{x}_i + \hat{s}_i) = 0$, $\kappa := \infty$ when $\sigma_m(A) = 0$, and $\|B^{-1}\| := \infty$ when $B$ lacks full rank.

We then establish the following upper bound for $\kappa \Phi$:
\begin{lemma}\label{lm:upperbound_kappa_phi_what}
	For the random LP defined in Definition \ref{def:random_lp}, if $\hat{x}$ and $\hat{s}$ satisfy the strict complementary slackness condition, the following inequality holds:
	\begin{equation}\label{eqoflm:upperbound_kappa_phi_what}
		\kappa \Phi \le \kappa\cdot \|B^{-1}\|\|A\| \cdot \varphi \ .
	\end{equation}
\end{lemma}
\begin{proof}{Proof.}
	When $B$ is full-rank, Lemma \ref{lm:full-rank-unique-optima} implies that Condition \ref{cond:unique_optima} holds with $\calX^\star=\{\hat{x}\}$ and $\calS^\star=\{\hat{s}\}$. In this case, by Lemma \ref{lm:upperbound_phi}, \eqref{eqoflm:upperbound_kappa_phi_what} holds.   When $B$ is not full-rank, Lemma \ref{lm:full-rank-unique-optima} implies that Condition \ref{cond:unique_optima} does not hold and by definition, $\Phi = \infty$. Simultaneously, $\|B^{-1}\| = \infty$ so inequality \eqref{eqoflm:upperbound_kappa_phi_what} still holds in this case.
\hfill\Halmos\end{proof}

\vspace{8pt}
\textbf{Roadmap of the proof of Theorem \ref{thm:complexity-random}.}
The above Lemma \ref{lm:upperbound_kappa_phi_what} indicates that in order to analyze the tail behavior of $\kappa \Phi$, we can study the tails of its components: $\kappa$, $\|B^{-1}\|\|A\|$ and $\varphi$. Therefore, below we study them in Steps 1 to 3, after which we will combine these results to derive a tail bound for $\kappa\cdot \|B^{-1}\|\|A\| \cdot \varphi$ (Step 4) and complete the proof of Theorem \ref{thm:complexity-random} and Corollary \ref{cor:complexity-random} (Step 5). 

\vspace{8pt}
\textbf{Step 1. Tail bound of $\kappa$.}
Before showing the result, we define some constants that will be frequently used in the following analysis:
\begin{equation}\label{def:c_1c_2c_3}
	c_1:=  C_{rv2}\wedge 6, \quad c_2:= 20\sigma_A C_{rv1}, \quad c_3:= 2c_2^2 \ .
\end{equation}
Note that here $c_1,c_2$ and $c_3$ depend only on $\sigma_A$, the sub-Gaussian parameter of $A$, and are independent of the distribution of $(\hat{x}, \hat{s})$ and the values of $m$ and $n$, because $C_{rv1}$ and $C_{rv2}$ are the constants in Lemma \ref{lm:min_singular_value_rectangular} that depend only on $\sigma_A$. 
\begin{lemma}\label{lm:bound_of_kappa}
	The following inequality holds for $\kappa$:
	\begin{equation}\label{eqoflm:bound_of_kappa}
		\Pr\left[\kappa \ge 2c_2\cdot \frac{n}{d+1}\right] \le \left(\tfrac{1}{2}\right)^{d+1}+2e^{-c_1 n }
	\end{equation}
\end{lemma} 

\begin{proof}{Proof.}
	We apply Lemma \ref{lm:min_singular_value_rectangular} and Lemma \ref{lm:max_singular_value_rectangular} to bound the largest and smallest singular values of $A$. For any $\eps>0$, let $E_1$ denote the event $\sigma_m(A) > \eps (\sqrt{n} - \sqrt{m-1})$ and let $E_2$ denote the event $\sigma_1(A) < 10\sigma_A\sqrt{n}$.
	In the event $E_1 \cap E_2$, it holds that
	$
		\kappa = \frac{\sigma_1(A)}{\sigma_m(A)} < \frac{10\sigma_A\sqrt{n}}{\eps (\sqrt{n} - \sqrt{m-1})} = \frac{10\sigma_A}{\eps}\cdot \frac{\sqrt{n}(\sqrt{n} + \sqrt{m-1})}{d+1} \le \frac{20\sigma_A}{\eps}\cdot \frac{n}{d+1}.
	$
	and thus
	\begin{equation}\label{eq:bound_of_kappa_1}
		\begin{aligned}
			\Pr\left[\kappa \ge \frac{20\sigma_A}{\eps}\cdot \frac{n}{d+1}\right] \le 1- \Pr\left[   E_1 \cap E_2 \right] \le \Pr\left[ E_1^c \right] + \Pr\left[ E_2^c \right]  \le (C_{rv1} \varepsilon)^{d+1}+e^{-C_{rv2} n} + e^{-6n}
		\end{aligned}
	\end{equation}
	for any $\eps>0$. 
	Here the last inequality is due to  Lemmas \ref{lm:min_singular_value_rectangular} and \ref{lm:max_singular_value_rectangular}.
		Therefore, setting $\eps=\frac{1}{2C_{rv1}}$ yields the desired inequality \eqref{eqoflm:bound_of_kappa} for $c_1$ and $c_2$ defined in \eqref{def:c_1c_2c_3}.
\hfill\Halmos\end{proof}

\vspace{8pt}
\textbf{Step 2. Tail bound of $\|B^{-1}\|\|A\|$.} 
Next, we analyze the tail bound of $\|B^{-1}\|\|A\|$.

\begin{lemma}\label{lm:bound_of_BinvA}
	For any $t > 0$, the following inequality holds:
	\begin{equation}\label{eqoflm:bound_of_BinvA}
		\Pr\Big[\|B^{-1}\|\|A\| \ge
		t
		\Big] \le   \frac{c_2\sqrt{mn}}{t}+2e^{-c_1 m} \ .
	\end{equation}
\end{lemma}

\begin{proof}{Proof.}
	The proof follows a similar structure to that of Lemma \ref{lm:bound_of_kappa}. 
	We use Lemmas \ref{lm:min_singular_value_rectangular} and  \ref{lm:max_singular_value_rectangular} to bound the largest singular value of $A$ (equal to $\|A\|$) and the smallest singular value of $B$ (equal to $1/\|B^{-1}\|$).
	For any $\eps>0$, let $E_1$ denote the event  $\sigma_m(B) > \eps (\sqrt{m} - \sqrt{m-1})$ and let $E_2$ denote the event $\sigma_1(A) < 10\sigma_A\sqrt{n}$.
	In the event $E_1 \cap E_2$, it holds that
	$
		\|B^{-1}\|\|A\| = \frac{\sigma_1(A)}{\sigma_m(B)} < \frac{10\sigma_A\sqrt{n}}{\eps (\sqrt{m} - \sqrt{m-1})} = \frac{10\sigma_A}{\eps}\cdot  \sqrt{n}(\sqrt{m} + \sqrt{m-1})  \le \frac{20\sigma_A}{\eps}\cdot \sqrt{mn} .
	$
	and thus
	\begin{equation}\label{eq:bound_of_BinvA_1}
		\begin{aligned}
			\Pr\left[\|B^{-1}\|\|A\| \ge \frac{20\sigma_A}{\eps}\cdot \sqrt{mn}\right] & \le 1- \Pr\left[   E_1 \cap E_2 \right] \le \Pr\left[ E_1^c \right] + \Pr\left[ E_2^c \right] \\
		                                                                                  & \le C_{rv1} \varepsilon +e^{-C_{rv2} m} + e^{-6n} \le C_{rv1} \varepsilon +e^{-C_{rv2} m} + e^{-6m}
		\end{aligned}
	\end{equation} where the third inequality uses Lemmas \ref{lm:min_singular_value_rectangular} and \ref{lm:max_singular_value_rectangular}, and the last inequality is due to $m\le n$.
	Finally,  replacing $\eps$ with $\frac{20\sigma_A}{t} \cdot \sqrt{mn}$  and using the definitions of $c_1$ and $c_2$ from \eqref{def:c_1c_2c_3}, we conclude the inequality \eqref{eqoflm:bound_of_BinvA}. The result is valid for any $\varepsilon > 0$ and the corresponding $t > 0$.
\hfill\Halmos\end{proof}

\vspace{8pt}
\textbf{Step 3. Tail bound of $\varphi$.}
Before showing the result, we define two constants $c_4, c_5> 0$ that will be frequently used:
\begin{equation}\label{def:c_3c_4}
	c_4 :=  \mu_u + 2 \sigma_u \quad \text{ and } \quad c_5 := C_{rv2}\wedge 2 \ .
\end{equation} 
Recall that here $\mu_u$ and $\sigma_u$ denote the maxima of the means and the sub-Gaussian parameters of $u$'s components, respectively. 

\begin{lemma}\label{lm:bounds_of_u}
	For all $t  >0$, the following tail bound holds:
	\begin{equation}\label{eqoflm:bounds_of_u_2}
		\Pr\left[\varphi \ge t\right] \le \frac{c_4n^2}{t}   + e^{-2n} \ .
	\end{equation}
	Furthermore, the following probabilistic bound holds:
	\begin{equation}\label{eqoflm:bounds_of_u_1}
		\Pr\left[\min_{1\le i\le n}(\hat{x}_i+\hat{s}_i)\le c_4\right]\ge1-e^{-2n} \ .
	\end{equation}
\end{lemma}
\begin{proof}{Proof.}

	Recall that the vector $u := (u^1, u^2)$ equals $\hat{x} + \hat{s}$ by definition, and
	\begin{equation}\label{eq:u_is_what}
		\min_{1 \le i \le n} (\hat{x}_i + \hat{s}_i) = \min_{1 \le i \le n}u_i \ , \ \ \varphi=\frac{\|\hat{x}+\hat{s}\|_1}{\min_{1 \le i \le n} (\hat{x}_i + \hat{s}_i)} = \frac{\sum_{i=1}^n u_i}{\min_{1 \le i \le n}u_i} \ .
	\end{equation}
	In the remainder of the proof we will mainly work on $u$.

	We first derive the probabilistic upper bound of $\sum_{i=1}^n u_i$.
	Suppose that for each $u_i$, its sub-Gaussian parameter is $\sigma_i$. Using Lemma \ref{lm:hoeffding}, it holds for all $t>0$ that
	\begin{equation}\label{eq:lm:bounds_of_u_1}
		\Pr\left[\sum_{i=1}^n u_i \geq n\mu_u+t\right] \leq \Pr\left[\sum_{i=1}^n u_i \geq \sum_{i=1}^n \E[u_i]+t\right]
		\le \exp \left\{-\frac{t^2}{2 \sum_{i=1}^n \sigma_i^2}\right\}
		\le \exp \left\{-\frac{t^2}{2 n \sigma_u^2}\right\} \ .
	\end{equation}
	Here the first inequality is due to $ \E[u_i] \le \mu_u$ for all $i=1,2,\dots,n$, the second inequality uses Lemma \ref{lm:hoeffding}, and the last inequality is due to $\sigma_u\ge \sigma_i$ for all $i=1,2,\dots,n$.
	Taking $t=2n\sigma_u$ in \eqref{eq:lm:bounds_of_u_1} gives
	$\Pr[\sum_{i=1}^n u_i\ge nc_4]\le e^{-2n}$. Since
	$\min_{1\le i\le n}u_i\le n^{-1}\sum_{i=1}^n u_i$, this proves \eqref{eqoflm:bounds_of_u_1}.

	Next we analyze the probabilistic lower bound of $\min_{1\le i \le n}u_i$. Because the probability density of each $u_i$ is bounded by $1$, the union bound gives, for any $t>0$,
	\begin{equation}\label{eq:lm:bounds_of_u_2}
		\Pr\left[\min_{1 \le i \le n}u_i \le t\right]
		\le \sum_{i=1}^n\Pr[u_i\le t]\le nt \ .
	\end{equation}
	Let $E_1$ denote the event $\sum_{i=1}^n u_i \ge n\mu_u+2n\sigma_u$. As shown above, $\Pr[E_1] \le e^{-2n}$. For all $\delta > 0$ we use $E_{2,\delta}$ to denote the event $\min_{1 \le i \le n}u_i \le \frac{\delta}{n}$. By \eqref{eq:lm:bounds_of_u_2} we have $\Pr[E_{2,\delta}] \le  \delta$.
	Therefore, for all $\delta >0$, it holds that
	\begin{equation}\label{eq:lm:bounds_of_u_3}
		\Pr\left[
			\frac{\sum_{i=1}^n u_i}{\min_{1\le i \le n} u_i} \ge \frac{n^2(\mu_u + 2\sigma_u)}{\delta}
			\right]  \le
		\Pr\left[E_1 \text{ or } E_{2,\delta}\right] \le \Pr[E_1] + \Pr[E_{2,\delta}] \le e^{-2n} + \delta \ .
	\end{equation}
	Replacing \(\delta\) in \eqref{eq:lm:bounds_of_u_3} with \(\frac{c_4 n^2}{t}\) and noting that $\mu_u + 2\sigma_u = c_4$ and $\frac{\sum_{i=1}^n u_i}{\min_{1 \le i \le n}u_i} = \varphi$, we can conclude \eqref{eqoflm:bounds_of_u_2}. 
\hfill\Halmos\end{proof}

\vspace{8pt}
\textbf{Step 4. Tail bound of $\kappa\Phi$.}
With the tail bounds of $\kappa$ and $\|B^{-1}\|\|A\|$ provided in Lemmas \ref{lm:bound_of_kappa} and \ref{lm:bound_of_BinvA}, we can derive a tail bound on their product. Using this tail bound, we then have the following tail bound of $\kappa\Phi$.
 
\begin{lemma}\label{lm:tail_of_kappaPhi}
	For any $\delta \in (0,1]$, the following bound holds:
	\begin{equation}\label{eqoflm:tail_of_kappaPhi_1}
		\Pr\left[
		\kappa\Phi \ge \frac{\ln(3/\delta)}{\delta}\cdot 6c_3c_4 \cdot \frac{ n^{3.5}m^{0.5}}{d+1}
		\right]\le \delta + \left(\tfrac{1}{2}\right)^{d+1} + 6e^{-c_5 m}  \ .
	\end{equation}
\end{lemma}
\begin{proof}{Proof.}

	First of all, we prove the following tail bound of $\kappa\|B^{-1}\|\|A\|$. For any $t > 0$, we claim the following:
	\begin{equation}\label{eqoflm:bound_of_kappaBinvA}
		\Pr\left[\kappa\|B^{-1}\|\|A\| \ge
		t
		\right] \le   \frac{1}{t}\cdot  \frac{c_3 n^{1.5}m^{0.5}}{d+1} + \left(\tfrac{1}{2}\right)^{d+1}+4e^{-c_1 m}  \ .
	\end{equation}
	Let $t_0$ be an arbitrary positive scalar, and we use $E_1$ and $E_2$ to denote the events $\kappa \ge \alpha_1:=  2c_2\cdot \frac{n}{d+1}$ and $\|B^{-1}\|\|A\| \ge \alpha_2:= t_0 \cdot c_2 \sqrt{mn}$, respectively.
	Then we have
	\begin{equation}\label{eq:lm:bound_of_kappaBinvA1}
		\begin{aligned}
			  & \Pr\left[\kappa\|B^{-1}\|\|A\| \ge t_0\cdot 2 c_2^2  \cdot \frac{n^{1.5}m^{0.5}}{d+1}  \right]
			=  \Pr\left[\kappa\|B^{-1}\|\|A\| \ge \alpha_1\alpha_2\right] \\
			\le \ & \Pr\left[\kappa \ge \alpha_1\text{ or }\|B^{-1}\|\|A\|\ge \alpha_2\right] = \Pr[E_1\cup E_2]\le \Pr[E_1] + \Pr[E_2] \\
			 \le\ &  \left(\tfrac{1}{2}\right)^{d+1}+2e^{-c_1 n}  +  \frac{1}{t_0}+2e^{-c_1 m}  \le \frac{1}{t_0} + \left(\tfrac{1}{2}\right)^{d+1}+4e^{-c_1 m} 
		\end{aligned}
	\end{equation}
	where the second inequality uses the union bound, and the third inequality uses Lemmas \ref{lm:bound_of_kappa} and \ref{lm:bound_of_BinvA} on $\Pr[E_1]$ and $\Pr[E_2]$, and the last inequality is due to $m\le n$. The above inequality holds for any $t_0>0$.
	Note that in the above \eqref{eq:lm:bound_of_kappaBinvA1}, $2c_2^2$ is equal to $c_3$. Let $t_0 = t\cdot \left(c_3  \cdot \frac{n^{1.5}m^{0.5}}{d+1}\right)^{-1}$ and then \eqref{eq:lm:bound_of_kappaBinvA1} simplifies to  \eqref{eqoflm:bound_of_kappaBinvA}. This proves the tail bound claim \eqref{eqoflm:bound_of_kappaBinvA}.

	With the tail bound of $\kappa\|B^{-1}\|\|A\|$ and the tail bound of $\varphi$ from Lemma \ref{lm:bounds_of_u}, we can now derive a tail bound of $\kappa\|B^{-1}\|\|A\|\cdot\varphi$ using Lemma \ref{lm:tail_of_product}.
	In the rest of the proof we use $Z$ to denote $\kappa \|B^{-1}\|\|A\|$. 
	
	From the definition of the probabilistic model, $Z$ and $\varphi$ are independent.  Then we can use Lemma \ref{lm:tail_of_product} to conclude the high-probability bound of $Z\varphi$. Specifically, for any $\delta \in (0,1]$:
	\begin{equation}\label{eq:lm:tail_of_kappaPhi_1}
		\Pr\left[
		Z\varphi \ge \frac{\ln(3/\delta)}{\delta}\cdot 6c_4n^2\cdot \frac{c_3 n^{1.5}m^{0.5}}{d+1}
		\right]\le \delta + \left(\tfrac{1}{2}\right)^{d+1} + 4e^{-c_1 m} + 2e^{-2n} \le \delta + \left(\tfrac{1}{2}\right)^{d+1} + 6e^{-c_5 m}\ .
	\end{equation}
	Here the last inequality holds by noting that $c_5:=  C_{rv2} \wedge 2$ is no larger than either $2$ or the $c_1$ defined in \eqref{def:c_1c_2c_3}, and  $n\ge m$.

	Finally, Lemma \ref{lm:upperbound_kappa_phi_what} states that $\kappa\Phi \le Z\varphi$ when $(\hat{x},\hat{s})$ satisfies strict complementary slackness, which holds almost surely. Therefore, for any $T > 0$, $\Pr[\kappa\Phi \ge T] \le \Pr[Z\varphi \ge T] $. Substituting this relationship into \eqref{eq:lm:tail_of_kappaPhi_1} completes the proof.
\hfill\Halmos\end{proof}

\vspace{8pt}
\textbf{Step 5. Proof of Theorem \ref{thm:complexity-random} and Corollary \ref{cor:complexity-random}.} Finally, we finish the proof.
\begin{proof}{Proof of Theorem \ref{thm:complexity-random}.}
	Let $E_1$ denote the event $\kappa\Phi \le \alpha_1:=\frac{\ln(3/\delta)}{\delta}\cdot 6c_3c_4 \cdot \frac{ n^{3.5}m^{0.5}}{d+1}$. Using \eqref{eq:basis_identification_T} of Theorem \ref{thm:closed-form-complexity}, $T_{basis}\le \check{c}_1  \kappa \Phi \cdot \ln(\kappa \Phi)$ for an absolute constant $\check{c}_1>1$, and we have:
	$$
	\Pr[T_{basis}\le \check{c}_1\alpha_1\ln(\check{c}_1\alpha_1)] \ge \Pr[\check{c}_1\kappa\Phi\cdot \ln(\kappa\Phi)\le  \check{c}_1\alpha_1\ln(\alpha_1)] = \Pr[E_1]\ge 1 - \Big(\delta + \left(\tfrac{1}{2}\right)^{d+1} + 6e^{-c_5 m}\Big) 
	$$
	where the first inequality uses $\check{c}_1 >1$ and the last inequality follows from Lemma \ref{lm:tail_of_kappaPhi}. 
	This inequality is exactly \eqref{eq:optimal_basis_random} if letting $c_0$ be $c_4$, $C_1$ be $6\check{c}_1 c_3$ and $C_0$ be $c_5$. Furthermore, $C_1$ and $C_0$ depend only (and at most polynomially) on $\sigma_A$.

	Now, let $E_2$ denote the event $\|B^{-1}\|\|A\| \le \alpha_2:=c_2\sqrt{mn}/\delta$ and let $E_3$ denote the event $\min_{1\le i\le n}u_i\le c_4$. By the convention following \eqref{eq:def_varphi}, $E_2$ implies that $B$ is full-rank. Since $u>0$ almost surely, Lemma \ref{lm:full-rank-unique-optima} then gives $x^\star=\hat x$ and $s^\star=\hat s$ on $E_2$. Thus, \eqref{eq:local_linear_convergence_T} gives
	$T_{local} \le \check{c}_2 \|B^{-1}\|\|A\|\max\{0,\ln(\min_{1\le i\le n}u_i/\eps)\}$.
	On $E_2\cap E_3$, this implies
	\begin{equation*}
		T_{local}\le \check c_2\alpha_2\max\left\{0,\ln\left(\frac{c_4}{\eps}\right)\right\}.
	\end{equation*}
	Furthermore, Lemma \ref{lm:bound_of_BinvA}, \eqref{eqoflm:bounds_of_u_1}, and the union bound give
	\begin{equation}\label{eq:thm:complexity-random-2}
		\Pr[E_2\cap E_3]\ge 1-\delta-2e^{-c_1m}-e^{-2n}
		\ge 1-\delta-3e^{-c_5m},
	\end{equation}
	where the last inequality uses $c_5\le c_1$, $c_5\le2$, and $m\le n$.
	Therefore, the preceding bound and \eqref{eq:thm:complexity-random-2} yield \eqref{eq:local_convergence_random} with $c_0:=c_4$, $C_2:=\check c_2c_2$, and $C_0:=c_5$.
\hfill\Halmos\end{proof}

\begin{remark}\label{rmk:constants-subgaussian}
From the above proof, we may conclude that the constants $C_0,C_1,C_2$ in Theorem~\ref{thm:complexity-random} can be explicitly represented as
$C_0 := (C_{\mathrm{rv}2}\wedge 2),$ $C_1 := 4800\,\check c_1\,\sigma_A^2\,C_{\mathrm{rv}1}^2,$ and $C_2 := 20\,\check c_2\,\sigma_A\,C_{\mathrm{rv}1}.$
Here $(C_{\mathrm{rv}1},C_{\mathrm{rv}2})$ are the constants from Theorem 1.1 of \citet{rudelson2009smallest} (used in Lemma~\ref{lm:min_singular_value_rectangular}), and $\check c_1,\check c_2$ are the absolute constants in Theorem 4.1 of \citet{xiong2024accessiblenew} (used in Theorem \ref{thm:closed-form-complexity}). Here $(C_{\mathrm{rv}1},C_{\mathrm{rv}2})$ depend only (and at most polynomially) on the sub-Gaussian parameter $\sigma_A$ of $A$.
\end{remark}

\begin{remark}\label{rmk:two-stage-linear-convergence} 
		The main idea of the proof of \eqref{eq:optimal_basis_random} is to establish the high-probability upper bound for $\kappa\Phi$ in Lemma \ref{lm:tail_of_kappaPhi}. Theorem 3.1 of \citet{xiong2024accessiblenew} shows a global linear convergence iteration bound $O(\kappa\Phi\cdot\log(\tfrac{\kappa\Phi(\|x^\star\|+\|s^\star\|)}{\eps}))$ for reaching an $\eps$-optimal solution. Since this complexity bound also heavily relies on $\kappa\Phi$, it is possible to use Lemma \ref{lm:tail_of_kappaPhi} to derive a high-probability bound on the global linear convergence of rPDHG. The bound in \eqref{eq:local_convergence_random} is consistent with fast local linear convergence because it depends only logarithmically on $c_0/\eps$.
\end{remark}

\begin{proof}{Proof of Corollary \ref{cor:complexity-random}.}
	Fix $\delta\in(0,1]$, set $\delta_0:=\delta/22$, and apply Theorem \ref{thm:complexity-random} with $\delta_0$. The condition in the corollary implies
	$2\delta_0>\max\{e^{-C_0m},2^{-(d+1)}\}$. Therefore, the right-hand sides of \eqref{eq:optimal_basis_random} and \eqref{eq:local_convergence_random} are lower bounded by $1-15\delta_0$ and $1-7\delta_0$, respectively. By the union bound, both events occur with probability at least $1-22\delta_0=1-\delta$.
	For $\eps\in(0,e^{-1}]$, $\max\{0,\ln(c_0/\eps)\}\le\big(1+\max\{0,\ln c_0\}\big)\ln(1/\eps)$. The factor $1+\max\{0,\ln c_0\}$ is absorbed into the $O(\cdot)$ notation, and the fixed rescaling from $\delta$ to $\delta_0$ is absorbed into both order terms.
\hfill\Halmos\end{proof}

\section{Proof of Theorem \ref{thm:complexity-gaussian}}
\label{sec:refined_analysis_gaussian}

In this section, we prove Theorem \ref{thm:complexity-gaussian}.
Section \ref{subsec:useful_helper_lemmas_gaussian} introduces a few lemmas on random matrices, such as probabilistic lower bounds for the intermediate singular values of a random matrix and the Hanson-Wright inequality. Section \ref{subsec:proof_of_thm_complexity_gaussian} contains the detailed proof of Theorem \ref{thm:complexity-gaussian}.

\subsection{Useful helper lemmas on random matrices}\label{subsec:useful_helper_lemmas_gaussian}

In this subsection, we introduce some useful lemmas. The first lemma below shows that the $k$-th largest singular value of a random matrix in $\R^{n\times n}$ grows linearly as $\Omega\left(\tfrac{n+1-k}{\sqrt{n}}\right)$ with high probability for all $k=1,2,\dots,n$. 

\begin{lemma}\label{lm:distribution_of_sigmai}
	Let $W\in\R^{n\times n}$ be a random matrix as in Definition \ref{def:random matrix}. For every $\delta > 0$, we have
	\begin{equation} \label{eq_of_lm:distribution_of_sigmai}
		\Pr\Big[\sigma_k(W) \ge \frac{\delta (n-k+1)}{4C_{rv1}\sqrt{n}}\text{ for all }k=1,2,\dots,n\Big] \ge 1 - \delta - n e^{-C_{rv2} n} \ ,
	\end{equation}
	where $C_{rv1}$ and $C_{rv2}$ are constants from Lemma \ref{lm:min_singular_value_rectangular}.
\end{lemma}

The proof of the above lemma uses the min-max principle for singular values and Lemma \ref{lm:min_singular_value_rectangular} on submatrices of the random matrix. The proof is deferred to Appendix \ref{appsec:useful_helper_lemmas_gaussian}.

The second lemma is the Hanson-Wright inequality, which provides a tail bound on the quadratic form of a random vector.

\begin{lemma}[Hanson-Wright inequality, Theorem 6.2.1 of \citet{vershynin2018high}]\label{lm:Hanson-Wright}
	Let $v \in \mathbb{R}^n$ be a random vector with i.i.d. components, each following a sub-Gaussian distribution with mean zero, unit variance, and parameter $\sigma$.
	Let $M$ be a matrix in $\R^{n\times n}$. Then for all $t \geq 0$, we have
	\begin{equation}\label{eq:lm:Hanson-Wright}
		\Pr\left[\left|v^{\top} M v-\E\left[v^{\top} M v\right]\right| \geq t\right] \leq 2\cdot \exp \left[-C_{hw} \cdot \min \left(\frac{t^2}{\sigma^4\|M\|_F^2}, \frac{t}{\sigma^2\|M\|}\right)\right] \ ,
	\end{equation} 
	where $C_{hw}$ is a positive absolute constant.
\end{lemma}

The Hanson-Wright inequality directly yields a tail bound for $\|W^{-1}v\|^2$, where $v$ is a random vector and $W$ is a matrix in $\mathbb{R}^{n \times n}$ with a guaranteed lower bound on all its singular values. The proof is deferred to Appendix \ref{appsec:useful_helper_lemmas_gaussian}.

\begin{lemma}\label{lm:tail_of_Binvv} 
	Suppose $W \in \mathbb{R}^{n \times n}$ satisfies the condition that there exists a constant $c_0 > 0$ such that  for all $k=1,2,\dots,n$ it holds that  
	\begin{equation}\label{eqoflm:tail_of_Binvv_1}
		\sigma_k(W) \ge \frac{c_0(n-k+1)}{\sqrt{n}} \ .
	\end{equation}
	Let $v \in \mathbb{R}^n$ be a random vector with i.i.d. components that follow a zero-mean, unit-variance, sub-Gaussian distribution of parameter $\sigma$. Then for all $\gamma \ge 0$, the following inequality holds:
	\begin{equation}\label{eqoflm:tail_of_Binvv_2}
\Pr\left[\|W^{-1}v\|^2 \ge \frac{2n(1+\sigma^2\gamma)}{c_0^2}  \right]
\le 2\cdot \exp\!\left(-2C_{hw}\min\{\gamma^2,\gamma\}\right) \ .
\end{equation}
	Here $C_{hw}$ is the constant in Lemma \ref{lm:Hanson-Wright}.
\end{lemma}

\subsection{Proof of Theorem \ref{thm:complexity-gaussian}}\label{subsec:proof_of_thm_complexity_gaussian}

In the proof of Theorem \ref{thm:complexity-random}, we used Lemma \ref{lm:upperbound_phi} to derive an upper bound for $\kappa\Phi$. However, directly using the expression of $\Phi$ in Definition \ref{def_geometric_measure}, we can derive a tighter upper bound for $\kappa\Phi$ in the case of the probabilistic model with Gaussian constraint matrices. The proof of Theorem \ref{thm:complexity-gaussian} is organized into four steps.

\vspace{8pt}
\textbf{Roadmap of the proof of Theorem \ref{thm:complexity-gaussian}.}
In Step 1, we show that $\kappa\Phi$ can be bounded as $\kappa\Phi \le \kappa \cdot \varphi \cdot (Z_p \vee Z_d)$, where $Z_p$ and $Z_d$ are newly defined random variables. We have already established probabilistic upper bounds for $\kappa$ and $\varphi$ in Section \ref{sec:prove_thm_random}, so in Step 2, we derive probabilistic upper bounds for $Z_p \vee Z_d$.  In Step 3, we combine these results to obtain a new high-probability bound for $\kappa\Phi$. Finally, we complete the proof of Theorem \ref{thm:complexity-gaussian} and Corollary \ref{cor:complexity-gaussian} in Step 4. 

\vspace{8pt}
\textbf{Step 1. A new upper bound on $\kappa\Phi$.} 
We start with some useful properties of Gaussian matrices. Gaussian matrices are orthogonally invariant. Specifically, suppose that $A$ is an $m\times n$ Gaussian matrix where $m < n$. The null space of $A$ has the same distribution as the row space of another $d \times n$ Gaussian matrix $Q$ that depends on $A$.  This property yields useful features of the probabilistic model, such as symmetry between the primal and dual problems (see more in \citet{todd1991probabilistic}). In other words, for an instance of the probabilistic model with a Gaussian constraint matrix, the dual problem \eqref{pro:dual LP_s} is itself an instance of the probabilistic model with a Gaussian constraint matrix.

\begin{lemma}[Theorem 2.4 of \citet{todd1991probabilistic}]
	For an instance of the probabilistic model in Definition \ref{def:random_lp} with a Gaussian constraint matrix $A \in \mathbb{R}^{m \times n}$, the dual problem \eqref{pro:dual LP_s} is also an instance of the probabilistic model with a Gaussian constraint matrix $Q \in \R^{d \times n}$. % Furthermore, the $\hat{x}$ and $\hat{s}$ are interchanged in the dual problem.
\end{lemma}

According to Lemma \ref{lm:full-rank-unique-optima}, instances of the probabilistic model with Gaussian constraint matrices satisfy Condition \ref{cond:unique_optima} almost surely with $\calX^\star = \{\hat{x}\}$ and $\calS^\star = \{\hat{s}\}$, because Gaussian matrices are full-rank almost surely. Consequently, the optimal basis of the primal problem is almost surely $\{1, 2, \dots, m\}$, while the optimal basis of its dual problem \eqref{pro:dual LP_s} is almost surely $\{m+1, m+2, \dots, n\}$. Let $\Theta$ denote the primal optimal basis $\{1, 2, \dots, m\}$, and let $\bar{\Theta}$ denote its complement $\{m+1, m+2, \dots, n\}$. Then $A_{\Theta}$ is exactly $B$ and $A_{\bar{\Theta}}$ is $N$.
For the primal problem, the simplex tableau $A_{\Theta}^{-1}A_{\bar{\Theta}}$ at the optimal solution is $B^{-1}N$. For the dual problem, it is $Q_{\bar{\Theta}}^{-1}Q_{\Theta}$, where $Q_{\bar{\Theta}}^{-1}Q_{\Theta}$ is equal to $-(B^{-1}N)^\top$, as stated below:

\begin{lemma}[Lemma 3.6 of \citet{xiong2024accessiblenew}]\label{lm:formula_QinvQ}
	The matrix $Q_{\bar{\Theta}}^{-1} Q_{\Theta}$ is equal to $-(B^{-1}N)^\top$.
\end{lemma}

Using Lemma \ref{lm:formula_QinvQ}, the terms $\|(B^{-1}N)_{i, \cdot}\|^2$ in Definition \ref{def_geometric_measure} of $\Phi$ can be converted into $\|(Q_{\bar{\Theta}}^{-1} Q_{\Theta})_{\cdot, i}\|^2$. Given that Condition \ref{cond:unique_optima} holds almost surely with $\calX^\star = \{\hat{x}\}$ and $\calS^\star = \{\hat{s}\}$, we can rewrite $\Phi$ as follows:
\begin{equation}
	\begin{aligned}
		\Phi &  =  \big(\|\hat{x}+\hat{s}\|_1\big) \cdot \max\left\{
		\max_{1 \le j \le d} \frac{\sqrt{\|(B^{-1}N)_{\cdot,j}\|^2+1} }{\hat{s}_{j + m}}, \ \max_{1\le i \le m} \frac{\sqrt{\|(Q_{\bar{\Theta}}^{-1} Q_{\Theta})_{\cdot,i}\|^2+1} }{\hat{x}_{i}}
		\right\}                                                                \\
		& \le \underbrace{\frac{\|\hat{x}+\hat{s}\|_1}{\min_{1\le i \le n}(\hat{x}_i+\hat{s}_i)}}_{\text{Denoted by $\varphi$}} \cdot \max\Bigg\{\underbrace{\max_{1 \le j \le d} \sqrt{\|(B^{-1}N)_{\cdot,j}\|^2+1}}_{\text{Denoted by $Z_p$}}, \ \underbrace{\max_{1\le i \le m}   \sqrt{\|(Q_{\bar{\Theta}}^{-1} Q_{\Theta})_{\cdot,i}\|^2+1}}_{\text{Denoted by $Z_d$}}  \Bigg\} 
	\end{aligned}
\end{equation}
in which the first term of the product is denoted by $\varphi$, and the two terms in the bracket are denoted by $Z_p$ and $Z_d$, respectively. Therefore, almost surely, we have
\begin{equation}\label{eq:tight_upperbound_phi}
	\kappa \Phi \le \kappa \cdot \varphi \cdot (Z_p \vee Z_d) \ .
\end{equation}
Tail bounds for $\kappa$ and $\varphi$ have already been established in Lemma \ref{lm:bound_of_kappa} and Lemma \ref{lm:bounds_of_u}, respectively. Thus, the primary focus will be on deriving a tail bound for $Z_p \vee Z_d$.

\vspace{8pt}
\textbf{Step 2. Tail bound of $Z_p \vee Z_d$.}
Although $Z_p$ and $Z_d$ are not independent, they have a symmetric structure, and the matrices $B$ and $Q_{\bar{\Theta}}$
are both Gaussian random matrices, which facilitates our analysis. To simplify notation, we define the following parameters for any $\gamma \ge 0$: 
	\begin{equation}\label{eq:def of c6 and delta_gamma}
		c_6 := \sqrt{256\cdot C_{rv1}^2 + 1}
		\quad \text{and} \quad
		\delta_\gamma := 2n \cdot \exp\!\left(-2C_{hw}\min\{\gamma^2,\gamma\}\right)
		+ n\cdot e^{-C_{rv2} \cdot (m \wedge d)} \ .
	\end{equation}
	Here $c_6$ is a fixed absolute constant as $C_{rv1}$ depends only on $\sigma_A$, which is equal to $1$ for Gaussian matrices. Additionally, $\delta_\gamma$ decreases exponentially with $\gamma$, $m$, and $d$, so for sufficiently large $\gamma$, $m$, and $d$, $\delta_\gamma$ becomes negligible compared to $\delta$. 

\begin{lemma}\label{lm:bounds of Zp} 
	For all $\delta >0$ and $\gamma \ge 1$,  the following bound holds:
	\begin{equation}\label{eqoflm:bounds of Zp simplified}
		\Pr\left[Z_p \vee Z_d \ge  \frac{c_6\sqrt{n\gamma}}{\delta}  \right]
		\le \delta + \delta_\gamma 
	\end{equation}
	for $c_6$ and $\delta_\gamma$ defined in \eqref{eq:def of c6 and delta_gamma}.
\end{lemma}

\begin{proof}{Proof.}
	When $\delta > 1$, \eqref{eqoflm:bounds of Zp simplified} is trivial. Later we consider the case $\delta \in (0,1]$.

	We first study the upper bound for $Z_p$. The same reasoning will apply symmetrically to $Z_d$.
	For the matrix $B$, let $E_{\delta}$ denote the event $\sigma_k(B)\ge \frac{\delta}{4C_{rv1}}\cdot\frac{m-k+1}{\sqrt{m}}$
	for all $k=1,2,\dots,m$. By Lemma \ref{lm:distribution_of_sigmai}, we have
	\begin{equation}\label{eq:bounds of Zp_1}
		\Pr\left[E_{\delta} \right] \ge 1 - \delta - m e^{-C_{rv2} m} \ .
	\end{equation}
	Next, Lemma \ref{lm:tail_of_Binvv} implies that for each $j\in\{1,\dots,d\}$,
	\begin{equation}\label{eq:bounds of Zp_2}
		\Pr\left[ \|B^{-1}N_{\cdot,j}\|^2 \ge \frac{2m(1+\sigma^2\gamma)}{\tfrac{\delta^2}{16C_{rv1}^2}}
		\;\middle|\; E_{\delta} \right]
		\le 2\cdot \exp\!\left(-2C_{hw}\min\{\gamma^2,\gamma\}\right) \ .
	\end{equation}
	Here $\sigma=1$ because $N_{\cdot,j}$ is a standard Gaussian vector (unit variance entries).
	Using the union bound over all $j \in \{1,\dots,d\}$, it follows that:
	\begin{equation}\label{eq:bounds of Zp_3}
		\Pr\left[ \max_{1\le j \le d}\|B^{-1}N_{\cdot,j}\|^2 \ge \frac{32mC_{rv1}^2(1+\gamma)}{\delta^2}
		\;\middle|\; E_{\delta} \right]
		\le 2d\cdot \exp\!\left(-2C_{hw}\min\{\gamma^2,\gamma\}\right) \ .
	\end{equation}
	Combining \eqref{eq:bounds of Zp_1} and \eqref{eq:bounds of Zp_3}, we have:
	\begin{equation}\label{eq:bounds of Zp_4}
		\Pr\left[\max_{1\le j \le d} \|B^{-1}N_{\cdot,j}\|^2 \ge \frac{32mC_{rv1}^2(1+\gamma)}{\delta^2}  \right]
		\le 2d\cdot \exp\!\left(-2C_{hw}\min\{\gamma^2,\gamma\}\right) + \delta + m e^{-C_{rv2} m}\ .
	\end{equation}
	Similarly, for $\max_{1\le i \le m} \|(Q_{\bar{\Theta}}^{-1}Q_{ \Theta})_{\cdot,i}\|^2$, we have:
	\begin{equation}\label{eq:bounds of Zp_6}
		\Pr\left[\max_{1\le i \le m} \|(Q_{\bar{\Theta}}^{-1}Q_{ \Theta})_{\cdot,i}\|^2
		\ge  \frac{32dC_{rv1}^2(1+\gamma)}{\delta^2}   \right]
		\le 2m\cdot \exp\!\left(-2C_{hw}\min\{\gamma^2,\gamma\}\right) + \delta + d e^{-C_{rv2} d}\ .
	\end{equation}
	Combining \eqref{eq:bounds of Zp_4} and \eqref{eq:bounds of Zp_6}, and noting that $n = m+d$, we have:
	\begin{equation}\label{eq:bounds of Zp_7}
		\Pr\left[Z_p \vee Z_d \ge \sqrt{\frac{32nC_{rv1}^2(1+\gamma)}{\delta^2} + 1}  \right]
		\le 2n \cdot \exp\!\left(-2C_{hw}\min\{\gamma^2,\gamma\}\right) + 2\delta + n\cdot e^{-C_{rv2} \cdot (m \wedge d)}\ .
	\end{equation} 
	Replacing $\delta$ in \eqref{eq:bounds of Zp_7} with $\delta/2$ proves:
	\begin{equation}\label{eqoflm:bounds of Zp}
	\Pr\left[Z_p \vee Z_d \ge \sqrt{\frac{128n\cdot C_{rv1}^2(1+\gamma)}{\delta^2} + 1}  \right]
	\le \delta + 2n \cdot \exp\!\left(-2C_{hw}\min\{\gamma^2,\gamma\}\right)
	+ n\cdot e^{-C_{rv2} \cdot (m \wedge d)} \ .
	\end{equation}

	To show \eqref{eqoflm:bounds of Zp simplified}, note that for $\gamma \ge 1$ and $\delta \in (0,1]$,
	\[
	\sqrt{\frac{128n\cdot C_{rv1}^2(1+\gamma)}{\delta^2} + 1 } \le 
	\sqrt{\frac{256  C_{rv1}^2 \cdot n\gamma}{\delta^2} + 1 } \le \sqrt{\frac{(256  C_{rv1}^2 + 1) \cdot n\gamma}{\delta^2}} = \frac{c_6\sqrt{n\gamma}}{\delta}.
	\]
	Substituting this inequality into \eqref{eqoflm:bounds of Zp}   yields \eqref{eqoflm:bounds of Zp simplified}.
\hfill\Halmos
\end{proof}

\vspace{8pt}
\textbf{Step 3. Tail bound of $\kappa\Phi$.}
We now combine the tail bounds established in Lemma \ref{lm:bound_of_kappa}, Lemma \ref{lm:bounds_of_u}, and Lemma \ref{lm:bounds of Zp}  to derive a tail bound for $\kappa \cdot \varphi \cdot (Z_p \vee Z_d)$, which directly provides a probabilistic upper bound for $\kappa\Phi$. To simplify the notation, we define two absolute constants:
\begin{equation}\label{eq:def_of_c7_and_c8}
	c_7:=  C_{rv2} \wedge \ln(2) \quad \text{and}\quad c_8:=  \max\big\{1,\tfrac{1}{\sqrt{2C_{hw}}} \big\} \ .
\end{equation}
Here $c_7$ is an absolute constant because $C_{rv2}$ only depends on the sub-Gaussian parameter $\sigma_A$, which is always equal to $1$ for Gaussian matrices.

\begin{lemma}\label{lm:refined_bound_of_kappaPhi}
For all $\delta \in (0,1]$, the following inequality holds:
\begin{equation}\label{eqoflm:refined_bound_of_kappaPhi}
\Pr\left[\kappa\Phi  \ge  24 c_2c_4 c_6 c_8 \cdot\frac{n^{3.5}}{d+1} \cdot
\frac{\ln\left(\tfrac{6}{\delta}\right)\sqrt{\ln\left(\frac{4n}{\delta}\right)}}{\delta} \right]
\le \delta + (n+5)\left(\frac{1}{e^{c_7}}\right)^{m \wedge d} \ .
\end{equation}
\end{lemma}

\begin{proof}{Proof.}
First of all, we study the tail bound of $\kappa\cdot(Z_p \vee Z_d)$. According to Lemma \ref{lm:bound_of_kappa}, using the union bound, we obtain the following probabilistic upper bound of $\kappa\cdot (Z_p \vee Z_d)$:
\begin{equation}\label{eq:refined_bound_of_kappaPhi_1}
	\Pr\left[\kappa\cdot (Z_p \vee Z_d) \ge  \frac{c_6\sqrt{n\gamma}}{\delta}\cdot2c_2\cdot \frac{n}{d+1}  \right]  \le \delta + \delta_\gamma + \underbrace{\left(\tfrac{1}{2}\right)^{d+1}+2\left(\tfrac{1}{e^{c_1}}\right)^n}_{\text{Denoted by $\delta_0$}}
\end{equation}	
for all $\gamma \ge 1$ and $\delta > 0$. Replacing $ \frac{c_6\sqrt{n\gamma}}{\delta}\cdot2c_2\cdot \frac{n}{d+1}$ by $t$ yields the following tail bound: 
\begin{equation}\label{eq:refined_bound_of_kappaPhi_2}
	\Pr\left[\kappa\cdot (Z_p \vee Z_d) \ge  t \right]  \le \frac{2c_2 c_6\sqrt{n\gamma}\cdot \frac{n}{d+1}}{t}  + \delta_\gamma +  \delta_0 \ .
\end{equation}	

Next, observe that $\kappa \cdot (Z_p \vee Z_d)$ and $\varphi$ are independent, and $\varphi$ already has a tail bound in Lemma \ref{lm:bounds_of_u}. Let $X = \kappa \cdot (Z_p \vee Z_d)$ and $Y = \varphi$. Using Lemma \ref{lm:tail_of_product}, we can derive the probabilistic upper bound for $XY = \kappa \cdot \varphi \cdot (Z_p \vee Z_d)$:
\begin{equation}\label{eq:refined_bound_of_kappaPhi_3}
	\Pr\left[\kappa\cdot\varphi \cdot (Z_p \vee Z_d)  \ge  \frac{6\cdot 2c_2 c_6\sqrt{n\gamma} \frac{n}{d+1} \cdot c_4n^2 \cdot \ln(3/\delta)}{\delta} \right]  \le \delta  + \delta_\gamma +  \delta_0 + 2\cdot e^{-2n}
\end{equation}	
for all $\delta \in (0,1]$ and $\gamma \ge 1$.   According to \eqref{eq:tight_upperbound_phi}, $\kappa\Phi \le\kappa\cdot\varphi \cdot (Z_p \vee Z_d)$ so
\begin{equation}\label{eq:refined_bound_of_kappaPhi_4}
	\Pr\left[\kappa\Phi  \ge  12 c_2c_4 c_6\cdot \frac{n^{3.5}}{d+1}   \cdot\frac{ \sqrt{\gamma} \cdot \ln(3/\delta)}{\delta} \right]  \le \delta  + \delta_\gamma +  \delta_0 + 2\cdot e^{-2n} \ .
\end{equation}

Then we simplify the right-hand side of \eqref{eq:refined_bound_of_kappaPhi_4}. Substituting the values of $\delta_\gamma$ and $\delta_0$, we have:
\begin{equation}\label{eq:refined_bound_of_kappaPhi_5}
	\begin{aligned}
	\delta  + \delta_\gamma + & \delta_0 + 2\cdot e^{-2n}   = \delta  + 2n \cdot \exp(-2C_{hw} \min\{\gamma^2,\gamma\}) +  n\cdot \left(\tfrac{1}{e^{C_{rv2}}}\right)^{m \wedge d}  +  \left(\tfrac{1}{2}\right)^{d+1} \\ & + 2\left(\tfrac{1}{e^{C_{rv2} \wedge 6}}\right)^n + 2\cdot \left(\tfrac{1}{e^2}\right)^{n} \le \delta  + 2n \cdot \exp(-2C_{hw} \min\{\gamma^2,\gamma\})  + (n+5)\left(\tfrac{1}{e^{c_7}}\right)^{m \wedge d }
	\end{aligned} 
\end{equation}	
where the last inequality follows from $c_7 = C_{rv2} \wedge \ln(2)$, as defined in \eqref{eq:def_of_c7_and_c8}. We now choose:
\begin{equation}\label{eq:refined_bound_of_kappaPhi_7}
\gamma = \max\left\{1,\; \frac{1}{2C_{hw}}\cdot \ln\left(\frac{2n}{\delta}\right) \right\} \ .
\end{equation}
In this way, since $\gamma \ge 1$ we have $\min\{\gamma^2,\gamma\}=\gamma$, and hence
$2n \cdot \exp\!\left(-2C_{hw}\min\{\gamma^2,\gamma\}\right)
= 2n \cdot \exp(-2C_{hw}\gamma) \le \delta$. From \eqref{eq:refined_bound_of_kappaPhi_5} we conclude that $\delta  + \delta_\gamma +  \delta_0 + 2\cdot e^{-2n} $ is further upper bounded by $ 2\delta     + (n+5)\left(\frac{1}{e^{c_7}}\right)^{m \wedge d}$. Substituting the value of $\gamma$ into \eqref{eq:refined_bound_of_kappaPhi_4} yields
\begin{equation}\label{eq:refined_bound_of_kappaPhi_8}
\Pr\left[\kappa\Phi  \ge 12 c_2c_4 c_6 \cdot\frac{n^{3.5}}{d+1} \cdot
\frac{\ln\left(\tfrac{3}{\delta}\right)\cdot
\max\left\{1,\; \sqrt{\frac{1}{2C_{hw}}\cdot \ln\left(\frac{2n}{\delta}\right)} \right\}}{\delta} \right]
\le 2\delta + (n+5)\left(\frac{1}{e^{c_7}}\right)^{m \wedge d} \ .
\end{equation}
Finally, since
$\max\left\{1, \sqrt{\frac{1}{2C_{hw}} \cdot \ln\left(\frac{2n}{\delta}\right)}\right\}
\leq c_8 \sqrt{\ln\left(\frac{2n}{\delta}\right)}$
for $\delta \le  \frac{1}{2}$ (as per \eqref{eq:def_of_c7_and_c8}), replacing $\delta$ with $\frac{\delta}{2}$ in the above inequality yields the desired probabilistic bound \eqref{eqoflm:refined_bound_of_kappaPhi}.
\hfill\Halmos\end{proof}

\vspace{8pt}
\textbf{Step 4. Proof of Theorem \ref{thm:complexity-gaussian} and Corollary \ref{cor:complexity-gaussian}.}
Finally, we finish the proof.

\begin{proof}{Proof of Theorem \ref{thm:complexity-gaussian}.}
	Let $E$ denote the event $\kappa\Phi \le  \alpha := 24 c_2c_4 c_6 c_8  \cdot\frac{n^{3.5}}{d+1}\cdot\frac{ \ln\left(\tfrac{6}{\delta}\right) \sqrt{\ln\left(\frac{4n}{\delta}\right)}}{\delta}.$
	By \eqref{eq:basis_identification_T} of Theorem \ref{thm:closed-form-complexity}, $T_{basis} \le \check{c}_1  \kappa\Phi \cdot \ln(\kappa \Phi)$ for an absolute constant $\check{c}_1 > 1$, and thus we have:
	\begin{equation}\label{eq:thm:complexity-gaussian-1}
		\Pr[T_{basis}\le \check{c}_1\alpha \cdot \ln(\check{c}_1\alpha)] \ge \Pr[\check{c}_1\kappa\Phi\cdot\ln(\kappa\Phi)\le \check{c}_1\alpha \cdot \ln(\alpha)] \ge \Pr[E] \ge 1 - \delta     - (n+5)\left(\tfrac{1}{e^{c_7}}\right)^{ m \wedge d} \ .
	\end{equation}
	where the first inequality uses $\check{c}_1 > 1$ and the last inequality is due to Lemma \ref{lm:refined_bound_of_kappaPhi}. 
	The inequality \eqref{eq:thm:complexity-gaussian-1} yields \eqref{eq:optimal_basis_gaussian} by setting $c_0=c_4$, $C_0=c_7$, and $C_1=24\check{c}_1 c_2c_6c_8$. Here $c_2,c_6,c_7$, and $c_8$ are all absolute constants because $\sigma_A = 1$ for a Gaussian matrix $A$, so $C_0$ and $C_1$ are absolute constants as well.  This proves \eqref{eq:optimal_basis_gaussian} with $C_0 = c_7$.
 
	As for $T_{local}$, the proof is identical to that of Theorem \ref{thm:complexity-random} because a Gaussian matrix is also a sub-Gaussian matrix. We may take $C_0=c_7$ because $c_1=C_{rv2}\wedge6\ge c_7=C_{rv2}\wedge\ln(2)$ and
	$1-\delta-2e^{-c_1m}-e^{-2n}\ge1-\delta-3e^{-c_7m}$.
	Thus, $C_0$ and $C_2=\check c_2c_2$ are absolute because $\sigma_A=1$. This completes the proof of Theorem \ref{thm:complexity-gaussian}.
\hfill\Halmos\end{proof}

\begin{remark}\label{rmk:constants-gaussian}
From the above proof, $(C_0,C_1,C_2)$ may be explicitly given by $C_0 := (C_{\mathrm{rv}2}\wedge \ln 2),$ $C_1 := 24\,\check c_1\,(20\sigma_A C_{\mathrm{rv}1})\sqrt{256C_{\mathrm{rv}1}^2+1}\cdot
\max\Big\{1,(2C_{\mathrm{hw}})^{-1/2}\Big\},$ and $
C_2 := 20\,\check c_2\,\sigma_A\,C_{\mathrm{rv}1}.$
Here $(C_{\mathrm{rv}1},C_{\mathrm{rv}2})$ are the constants from Theorem 1.1 of \citet{rudelson2009smallest},
$\check c_1,\check c_2$ are the absolute constants in Theorem 4.1 in \citet{xiong2024accessiblenew}, and $C_{\mathrm{hw}}$ is the absolute constant in the Hanson-Wright inequality as stated, e.g., in Theorem~6.2.1 of \citet{vershynin2018high}
(used in Lemma~\ref{lm:Hanson-Wright}). Note that $(C_{\mathrm{rv}1},C_{\mathrm{rv}2})$ depend only on $\sigma_A$ but for the Gaussian matrix $A$ defined in Definition \ref{def:gaussian_matrix}, $\sigma_A $ is equal to $1$. Thus, $C_0,C_1$, and $C_2$ are absolute constants.
\end{remark}

\begin{proof}{Proof of Corollary \ref{cor:complexity-gaussian}.}
	Fix $\delta\in(0,1]$, set $\delta_0:=\delta/8$, and apply Theorem \ref{thm:complexity-gaussian} with $\delta_0$. The condition in the corollary, $m\wedge d\le m$, and $n+5\ge6$ imply $(n+5)e^{-C_0(m\wedge d)}<4\delta_0$ and $3e^{-C_0m}<2\delta_0$. Thus, the right-hand sides of \eqref{eq:optimal_basis_gaussian} and \eqref{eq:local_convergence_gaussian} are lower bounded by $1-5\delta_0$ and $1-3\delta_0$, respectively. By the union bound, both events occur with probability at least $1-8\delta_0=1-\delta$. For $\eps\in(0,e^{-1}]$, $\max\{0,\ln(c_0/\eps)\}\le\big(1+\max\{0,\ln c_0\}\big)\ln(1/\eps)$. The factor $1+\max\{0,\ln c_0\}$ is absorbed into the $O(\cdot)$ notation, and the fixed rescaling from $\delta$ to $\delta_0$ is absorbed into both order terms.
\hfill\Halmos\end{proof}

\section{Experimental Confirmation of the High-Probability Polynomial-Time Complexity}\label{sec:experiments}

This section presents the experimental results of rPDHG applied to randomly generated LP instances to validate our high-probability complexity analysis. Section \ref{subsec:dependence_on_delta} confirms the tail behavior of the iteration counts in both stages. Section \ref{subsec:dependence_on_n_practice} demonstrates the polynomial dependence of the number of iterations on the dimension $n$.

We implement rPDHG for the standard-form problem \eqref{pro:primal LP}, following Algorithm \ref{alg: PDHG with restarts} in Appendix \ref{appsec:restartPDHG}. In our numerical experiments, we regard an LP instance as successfully solved if rPDHG computes a primal-dual solution pair $(x^k,y^k)$ within Euclidean distance $10^{-4}$ of the optimal solution, namely, $\|(x^k,y^k)-(x^\star,y^\star)\|\le 10^{-4}$.  Compared with the commonly used KKT residuals, Euclidean distance is a more straightforward measure of optimality, and the probabilistic model with Gaussian input data ensures easy access to $(x^\star,y^\star)$. 
It is also the same criterion used in the experiments of \citet{xiong2024accessiblenew}.
We manually classify the iterations into two stages. Stage I comprises all iterations until the support set of $x^k$ settles down to the optimal basis, while Stage II consists of all subsequent iterations. These stages correspond to the optimal basis identification and fast local convergence stages analyzed in Section \ref{sec:high-probability}.

\subsection{Tail behavior of the iteration count}\label{subsec:dependence_on_delta}

\begin{figure}[htbp]
		\FIGURE
		{
		\hspace*{\fill} 
		\subcaptionbox{\small Stage I \vspace{-5pt} \label{fig:stage1}}
		{\includegraphics[width=0.47\linewidth]{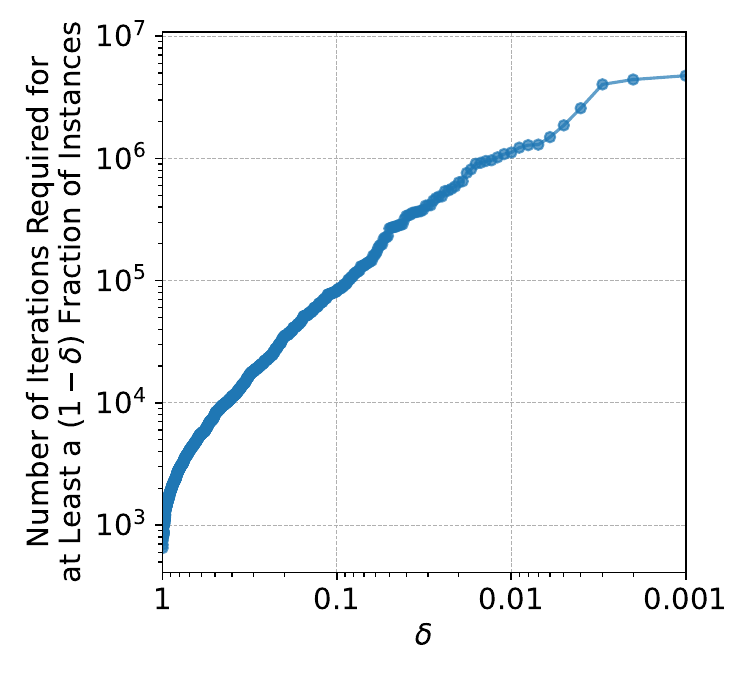}}
		\subcaptionbox{\small Stage II \vspace{-5pt}
		\label{fig:stage2}}
			{\includegraphics[width=0.47\linewidth]{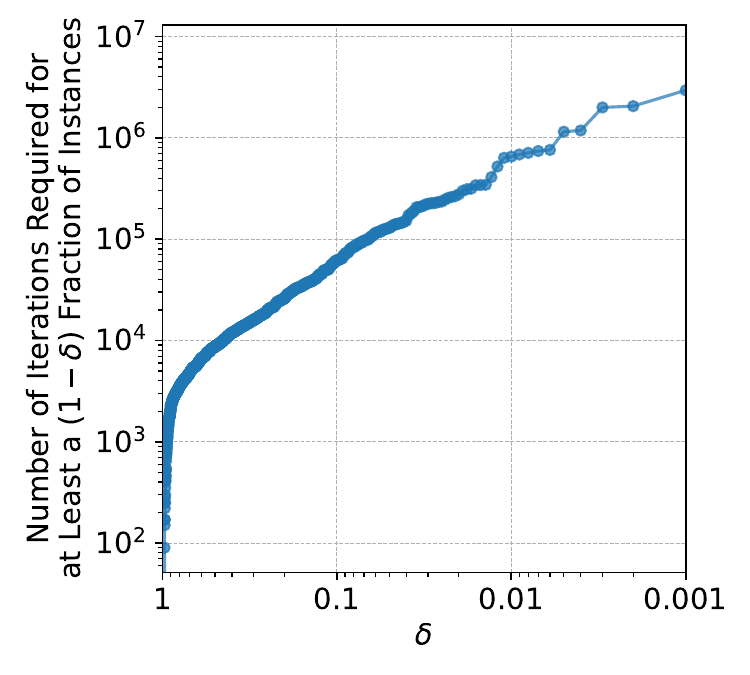}}
		\hspace*{\fill} % Add blank space on the right
		}
		{\normalsize Number of iterations required to complete Stages I and II of at least a $(1-\delta)$ fraction of instances \vspace{5pt} \label{fig:stages}}
		{} 
	\end{figure}

This subsection presents experimental confirmation of the tail behavior of the Stage-I and Stage-II iteration counts. We apply rPDHG to 1,000 LP instances generated according to Definition \ref{def:random_lp} for $n = 100$ and $m = 50$. The constraint matrix $A$ is a Gaussian matrix.
The vectors $\hat{x}$ and $\hat{s}$ are constructed with i.i.d. nonzero components, each drawn from the folded Gaussian distribution (absolute value of a Gaussian random variable of zero mean and unit variance). In all experiments, we use the objective vector $\bar c$ so that $A\bar c=0$.
Figure \ref{fig:stages} shows the number of iterations required to complete Stages I and II of at least a $(1-\delta)$ fraction of instances, for different values of $\delta$ in $(0,1)$.
	
Theorems \ref{thm:complexity-random} 
and \ref{thm:complexity-gaussian} theoretically predict that the number of iterations required to complete Stages I and II with a $(1-\delta)$ success rate should grow with $\frac{1}{\delta}$ (at rates $\widetilde{O}(\frac{1}{\delta})$ and $O(\frac{1}{\delta})$, respectively), except for exponentially small $\delta$. 
If the theory is exact, in Figures \ref{fig:stage1} and \ref{fig:stage2} (log-log plots with reversed horizontal axes), the slopes should be roughly equal to $1$, meaning that the number of iterations required has a reciprocal relationship with $\delta$. Indeed, Figure \ref{fig:stages} shows that both slopes are approximately $1$, particularly for $\delta$ between $0.1$ and $0.01$. These observations confirm our theoretical prediction for the dependence on $\delta$ and the tail behavior of the iteration counts.

\subsection{Polynomial dependence of the iteration count on the problem dimension}\label{subsec:dependence_on_n_practice} 
	
\begin{figure}[htbp]
		\FIGURE
		{
		\hspace*{\fill} 
		\subcaptionbox{\small Stage I \vspace{-5pt} \label{fig:dimensiontest_stage1}}
		{\includegraphics[width=0.4\textwidth]{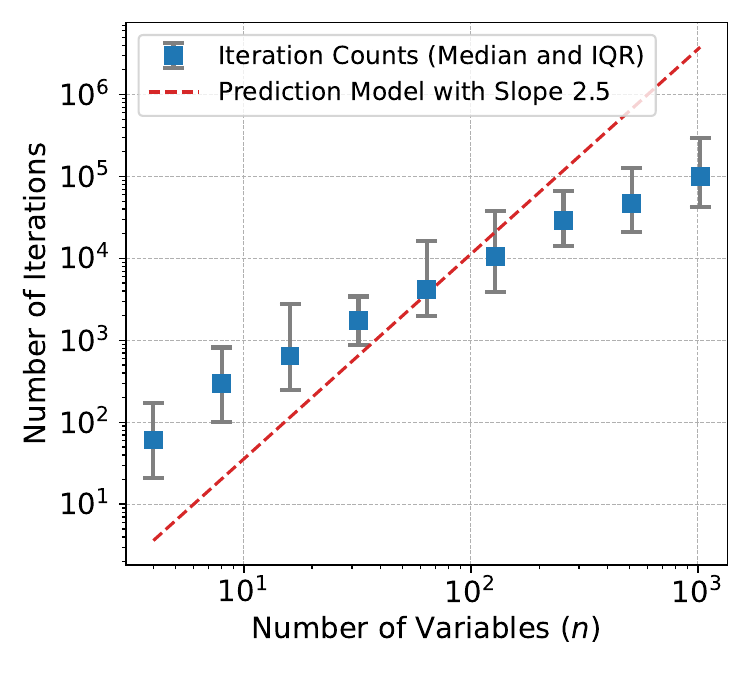}}
		\subcaptionbox{\small Stage II \vspace{-5pt}
		\label{fig:dimensiontest_stage2}}
			{\includegraphics[width=0.4\textwidth]{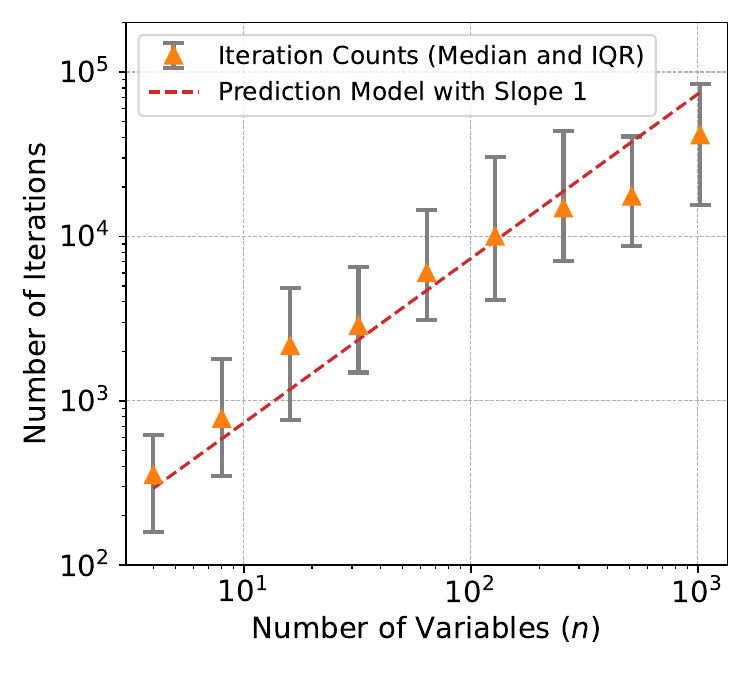}}
		\hspace*{\fill} % Add blank space on the right
		}
		{\normalsize The median values and interquartile range (IQR) of Stage-I and Stage-II iteration counts for various dimensions $n$ (for LP instances with Gaussian input data). \vspace{5pt} \label{fig:dependence_on_n}}
		{} 
	\end{figure}

This subsection examines how iteration counts of the two stages scale with the problem dimension (number of variables) in practice. We generate LP instances according to Definition \ref{def:random_lp} with $m=n/2$ for $n$ in $\{4,8,16,32,\dots\}$. As in Section \ref{subsec:dependence_on_delta}, the constraint matrix $A$ is a Gaussian matrix, and the vector $u$ has i.i.d. components, each obeying the folded Gaussian distribution.
For each value of $n$, we generate $100$ LP instances and compute the first quartile, median, and third quartile of both Stage-I and Stage-II iteration counts. Figure \ref{fig:dependence_on_n} shows the relation between these statistics of the iteration counts and the number of variables $n$.

Theorem \ref{thm:complexity-gaussian} predicts that in most cases, the Stage-I and Stage-II iteration counts should scale as $\widetilde{O}(n^{2.5})$ and $O(n)$, respectively.
If the theory is exact, since Figures \ref{fig:dimensiontest_stage1} and \ref{fig:dimensiontest_stage2} are log-log plots, the slopes should be approximately $2.5$ and $1$, respectively, corresponding to polynomial dependencies of degree $2.5$ and $1$ on $n$.
Based on these insights, we fit prediction models to the median iteration counts using slopes of $2.5$ and $1$ for Stages I and II. The predicted iteration counts using these models are shown by the dashed lines in the two figures for Stage I and Stage II, respectively. 
Figure \ref{fig:dimensiontest_stage1} shows that the Stage-I iteration count grows more slowly than the theoretical bound $\widetilde{O}(n^{2.5})$ in Theorem \ref{thm:complexity-gaussian}, as it does not grow as fast as the dashed line. Instead, the practical iteration count dependence on $n$ is smaller than our model by roughly a factor proportional to $n$. 
For Stage II, Figure \ref{fig:dimensiontest_stage2} shows that the iteration count aligns well with the $O(n)$ dependence predicted by Theorems \ref{thm:complexity-random} and \ref{thm:complexity-gaussian}. Overall, the slower growth of Stage-II iterations relative to Stage-I iterations validates the insight that local linear convergence is faster, particularly for larger problem dimensions.

A similar gap between observed practical performance and probabilistic analysis is common in the average-case analysis literature on classic methods as well. For instance, \citet{shamir1987efficiency} summarizes many reported observations of simplex method performance on real-life LP problems and concludes that the number of iterations for real-life problems is usually observed to be no higher than $O(m)$ for LP instances with $m$ linear constraints, but the average-case analyses of \citet{todd1986polynomial,adler1985simplex} and \citet{adler1986family} all prove a quadratic dependence of the expected iteration count on the problem dimension in theory. Similarly, the interior-point method iteration count often shows logarithmic growth in the problem dimension $n$ in practice, despite the polynomial dependence on $n$ proven in the probabilistic analyses by \citet{ye1994toward,anstreicher1999probabilistic,dunagan2011smoothed} and others. Our high-probability Stage-I iteration bound (shown by Theorem \ref{thm:complexity-gaussian}) also overestimates the empirical scaling by a factor of $O(n)$, while our Stage-II bound accurately predicts the observed dependence on $n$. In addition to revealing that the complexity of rPDHG is polynomial in dimension with high probability, an interesting research question is how to further shrink the gap between theory and practice in terms of the polynomial dependence on $n$.
We conjecture that the Stage-I iteration count might actually be bounded by $\widetilde{O}(\frac{n^{1.5}}{\delta})$ with probability at least $1-\delta$ for $\delta$ that is not exponentially small. Clearly, this issue merits more research in the future.

From Figure \ref{fig:dimensiontest_stage1}, a small prefactor (about $1.906$) multiplying $n^{2.5}$ already upper-bounds the median Stage-I iteration counts over the tested range. Similarly, Figure \ref{fig:dimensiontest_stage2} suggests that a modest prefactor (about $135.0$) multiplying $n$ already upper-bounds the median Stage-II iteration counts. These observations empirically suggest that the constants $c_0C_1$ and $C_2$ need not be large for the bounds in Theorem \ref{thm:complexity-gaussian} to hold with probability at least one half.

\begin{figure}[htbp]
		\FIGURE
		{
		\hspace*{\fill} 
		\subcaptionbox{\small Stage I \vspace{-5pt} \label{fig:dimensiontest_stage1_subgaussian}}
		{\includegraphics[width=0.4\textwidth]{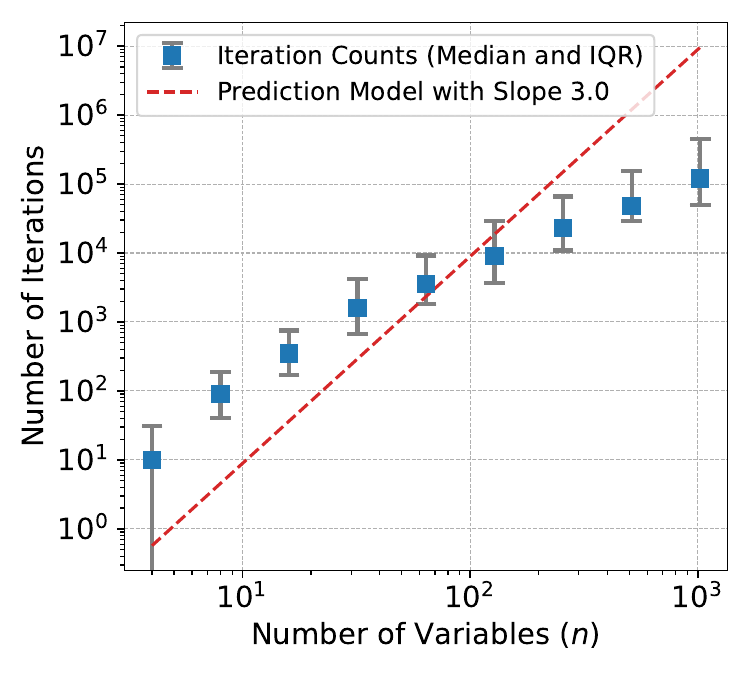}}
		\subcaptionbox{\small Stage II \vspace{-5pt}
		\label{fig:dimensiontest_stage2_subgaussian}}
			{\includegraphics[width=0.4\textwidth]{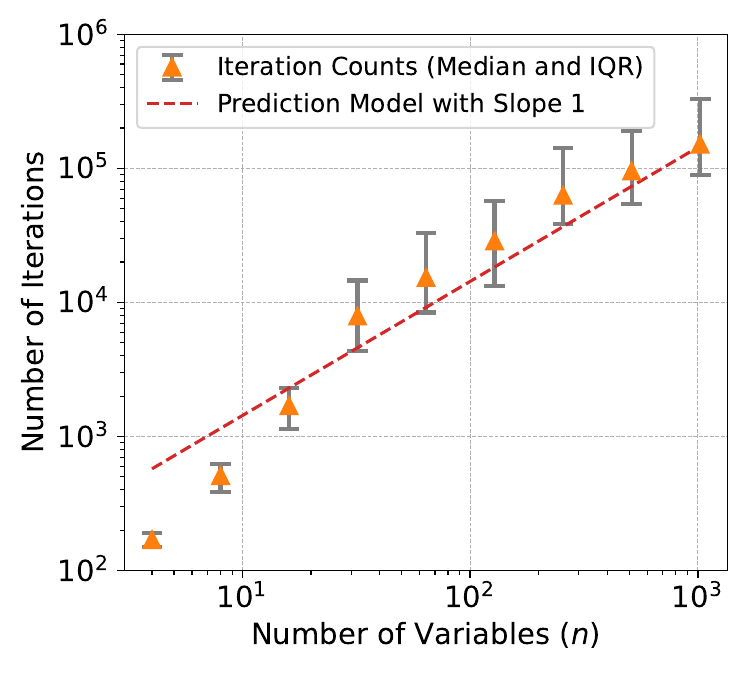}}
		\hspace*{\fill} % Add blank space on the right
		}
		{\normalsize The median values and interquartile range (IQR) of Stage-I and Stage-II iteration counts for various dimensions $n$ (for LP instances with sub-Gaussian input data). \vspace{5pt} \label{fig:dependence_on_n_sub-Gaussian}}
		{} 
	\end{figure}

Furthermore, we conducted the same experiment for LP instances with sub-Gaussian input data. The only difference is that the entries of 
$A$ are i.i.d. Rademacher random variables (taking values $1$ and $-1$ with probability $1/2$). The experimental results are shown in Figure \ref{fig:dependence_on_n_sub-Gaussian}. Theorem \ref{thm:complexity-random} predicts that in most cases, the Stage-I and Stage-II iteration counts should scale as $\widetilde{O}(n^{3})$ and $O(n)$, respectively. 
Based on these insights, we also fit prediction models to the median iteration counts using slopes of $3$ and $1$ for Stages I and II. The insights from this sub-Gaussian input data case (including the gap with the theoretical predictions) are no different from those in the experiment with Gaussian input data.

\section{How the Disparity among Optimal Solution Components Affects the Performance of rPDHG}\label{sec:effect_of_optimal_solution}

In this section, we use probabilistic analysis to investigate how the disparity among the optimal solution components affects the performance of rPDHG. Let $E_{\mathrm{uniq}}$ denote the event that Condition \ref{cond:unique_optima} holds, namely, that the primal and dual problems have unique optimal solutions. % Since each component of the nonnegative vector $u$ is absolutely continuous, $u>0$ almost surely. Thus, Lemmas \ref{lm:full-rank-unique-optima} and \ref{lm:full-rank-probability} imply $\Pr[E_{\mathrm{uniq}}] \ge 1-e^{-\nu m}>0$, so all conditional probabilities in this section are well-defined. 
On $E_{\mathrm{uniq}}$, let $x^\star$ and $s^\star$ denote the unique primal optimal solution and dual optimal slack, respectively, and define their disparity ratio as follows:
\begin{equation}\label{eq:def_of_phi}
	\phi := \frac{\frac{1}{n}\cdot\sum_{i=1}^n \left(x_i^\star + s_i^\star\right)}{\min_{1 \le i \le n} (x_i^\star + s_i^\star)}.
\end{equation}
Note that $\phi \geq 1$, with equality holding if and only if $x^\star + s^\star$ is a scalar multiple of the all-ones vector $e$.
The quantities in the conditional bounds below are evaluated on $E_{\mathrm{uniq}}$. Their values outside this event are immaterial. We conduct probabilistic analysis of the iteration bound conditioned on $E_{\mathrm{uniq}}$. Our result shows that, among instances with unique primal and dual optima, $T_{basis}$ depends linearly on the realized disparity ratio $\phi$ with high conditional probability, while $T_{local}$ is nearly independent of $\phi$.

\begin{theorem}\label{thm:complexity-conditioned}
	Suppose that rPDHG is applied to an instance of the probabilistic model in Definition \ref{def:random_lp} (with objective vector $\bar{c}$).  There exist constants $C_0,C_1,C_2>0$ that depend only (and at most polynomially) on $\sigma_A$ for which the following high-probability iteration bounds hold:
	\begin{enumerate}
		\item\text{(Optimal basis identification)} 
		\begin{equation}\label{eq:optimal_basis_conditioned}
			      \Pr\left[ 
			      T_{basis}
			      \le \frac{m^{0.5} n^{2.5}}{d+1}\cdot
			      \frac{C_1 \phi}{\delta}\cdot \ln\left(
					\tfrac{m^{0.5} n^{2.5}}{d+1}\cdot
					\tfrac{C_1 \phi }{\delta} 
			      \right)
			      \;\middle|\; E_{\mathrm{uniq}}
			      \right]\ge 1 - \delta - 4\left(\tfrac{1}{e^{C_0}}\right)^m - \left(\tfrac{1}{2}\right)^{d+1}
		      \end{equation}
		      for any $\delta \in (0,1]$.  
		\item\text{(Fast local convergence)} Let $\eps > 0$ be any given tolerance.  
		\begin{equation}\label{eq:local_convergence_conditioned}
			      \Pr\left[ 
			      T_{local}
			      \le m^{0.5} n^{0.5}
				  \cdot
			      \frac{C_2  }{\delta} 
			      \cdot \max\left\{0,\ln\left(
			      \tfrac{\min_{1 \le i \le n}(x_i^\star + s^\star_i)}{\eps}
			      \right)\right\}
			      \;\middle|\; E_{\mathrm{uniq}}
			      \right]\ge 1 - \delta - 2\left(\tfrac{1}{e^{C_0}}\right)^m
		      \end{equation}
		      for any $\delta \in (0,1]$.
	\end{enumerate} 
\end{theorem}

The above theorem characterizes the high-probability performance of rPDHG conditioned on the primal and dual optima being unique. The quantities $\phi$ and $\min_i(x_i^\star+s_i^\star)$ are their realized values on $E_{\mathrm{uniq}}$. 
Compared with Theorem \ref{thm:complexity-random}, the constant $c_0$ is absent, and the disparity ratio $\phi$ is the only quantity related to the optimal solution in \eqref{eq:optimal_basis_conditioned}. 

The quantity $\frac{n^{2.5}}{d+1}$ is at most as large as $\frac{n^{2.5}}{2}$.  When $d$ (recall $d = n-m$) is not too small compared to $n$ in the sense that $d + 1 \ge \frac{n}{C_3}$ for some absolute constant $C_3 >1$, the quantity $\frac{n^{2.5}}{d+1}$ is $O(n^{1.5})$. 
The high-probability iteration bound for Stage I (optimal basis identification) depends linearly on $\phi$, while that for Stage II (local convergence) depends only logarithmically on $\min_{1 \le i \le n}(x^\star_i + s^\star_i)$.  Roughly speaking, the greater the disparity among the components of $x^\star + s^\star$, the larger $\phi$ becomes, making an LP instance more challenging for Stage I of rPDHG according to our theory. For example, when $x^\star + s^\star = e$ (the all-ones vector) and $m$ is not too close to $n$, the Stage-I iteration count is bounded by $\widetilde{O}(\frac{n^{1.5} m^{0.5}}{\delta})$ with probability at least $1 - O(\delta)$ for non-exponentially small $\delta$. This bound on $T_{basis}$ is almost $n$ times smaller than the bound of $T_{basis}$ for the random LP instances in Theorem \ref{thm:complexity-random}. Conversely, highly imbalanced components in $x^\star + s^\star$ can lead to larger bounds on $T_{basis}$ than those in Theorem \ref{thm:complexity-random}.

It should be noted that metrics similar to the disparity defined above (such as the ratios between the largest and smallest components, or the smallest positive component) also influence the complexity of other methods, especially interior-point methods, for solving LPs. The proximity measures $\frac{\max_i x_i s_i}{\min_i x_i s_i}$ quantify the distance from the central path and enter the design and analysis of interior-point methods (see \citet{terlaky2009twenty,guler1993convergence}). The minimum positive coordinate scale over strictly complementary solutions is used to determine when the iterates of interior-point methods can identify the optimal face (see, e.g., \citet{terlaky2009twenty,ye1992finite}). A quantity involving the ratio between the norm and the smallest positive component of a pair of primal-dual strictly complementary optimal solutions directly appears in the finite termination analysis of interior-point methods (see, e.g., \citet{anstreicher1999probabilistic}). For these methods, such metrics typically enter the iteration complexity only logarithmically, unlike that for rPDHG.
Beyond interior-point methods, the strongly polynomial complexity analyses of the simplex method and policy iteration methods for discounted Markov decision problems rely on boundedness of ratios such as $\frac{\|x^\star\|_1}{\min_{i:\,x_i^\star>0} x_i^\star}$ (see, e.g., \citet{ye2011simplex}). Additionally, \citet{lu2024geometrynew} demonstrate that PDHG (without restarts) exhibits faster local linear convergence within a neighborhood whose size relates to $\min_{1 \leq i \leq n}\{x_i^{\star}+s_i^{\star}\}$. In these two settings, such metrics have a larger effect because they appear outside a logarithm.

\subsection{Application: Generating difficult LP instances}\label{subsec:depend_on_varphi}

Our analysis suggests that we can generate LP instances with controlled difficulty levels for rPDHG by manipulating the distributions of $\hat{x}$ and $\hat{s}$ in our probabilistic model. By Lemma \ref{lm:full-rank-unique-optima}, $\hat{x}$ and $\hat{s}$ are almost surely the unique primal and dual optimal solutions if the constraint matrix is a Gaussian matrix. The larger the ratio $\frac{\frac{1}{n}\sum_{i=1}^n (\hat{x}_i + \hat{s}_i)}{\min_{1 \le i \le n} (\hat{x}_i + \hat{s}_i)}$, the more difficult the LP instance is likely to be. Such instances can serve as useful benchmarks for evaluating the performance of rPDHG and other first-order LP methods.  Below we present an approach to generate difficult random LP instances.

Following Definition \ref{def:random_lp}, we consider random instances with a Gaussian constraint matrix $A \in \mathbb{R}^{m \times n}$, where $n=2m$.  The only difference from Definition \ref{def:random_lp} is that we define $\hat{x}=\hat{x}^l$ and $\hat{s}=\hat{s}^l$ for $l\ge0$ to control $\phi$: 
\begin{equation}
    \hat{x}^l:=
    \begin{pmatrix}
        u^l\\
        0_m
    \end{pmatrix}
    \quad \text{ and }\quad 
    \hat{s}^l:=
    \begin{pmatrix}
        0_m\\
        u^l
    \end{pmatrix} 
    \quad \text{ where }
    u^l:=\big[\underbrace{4^{-l},4^{-l},\ldots,4^{-l}}_{\lfloor m/2 \rfloor \text{ copies}},
    \underbrace{1, \ 1,\ \ldots,\ 1}_{m-\lfloor m/2 \rfloor \text{ copies}}\big] \ .
\end{equation}
Here $\lfloor m/2 \rfloor$ denotes the largest integer no greater than $m/2$.
These instances have full row rank almost surely, and $\hat{x}^l$ and $\hat{s}^l$ satisfy strict complementary slackness, ensuring that they are the unique pair of optimal solutions. 
Therefore, the smallest component of $x^\star+s^\star$ is 
\begin{equation}\label{def:smallestnnz_l}
	\min_{1 \le i \le n} (x_i^\star + s_i^\star)  = \min_{1 \le i \le n} (\hat{x}_i^l + \hat{s}_i^l)  = 4^{-l} \ ,
\end{equation}
and the corresponding disparity ratio $\phi_l$ is equal to
\begin{equation}\label{def:phi_l}
	\phi_l = \frac{1}{2m} \cdot \frac{2 \lfloor m/2 \rfloor \cdot 4^{-l} + 2 (m-\lfloor m/2 \rfloor)\cdot 1}{4^{-l}} = \frac{\lfloor m/2 \rfloor}{m} + \left(1-\frac{\lfloor m/2 \rfloor}{m}\right)\cdot 4^l \ .
\end{equation} 
Note that $\phi_l$ is approximately equal to $2^{2l-1}$ when $l$ and $m$ are large enough. In this way, the parameter $l$ controls the magnitude of $\phi_l$. Instances with large values of $l$ have large $\phi_l$, unbalanced $x^\star + s^\star$ and small $\min_{1 \le i \le n} (x_i^\star + s_i^\star)$.

\begin{figure}[htbp] 
		\FIGURE
		{
			\hspace*{\fill} 
		\subcaptionbox{\small Stage I \vspace{-5pt} \label{fig:phi_stage1}}
		{\includegraphics[width=0.47\textwidth]{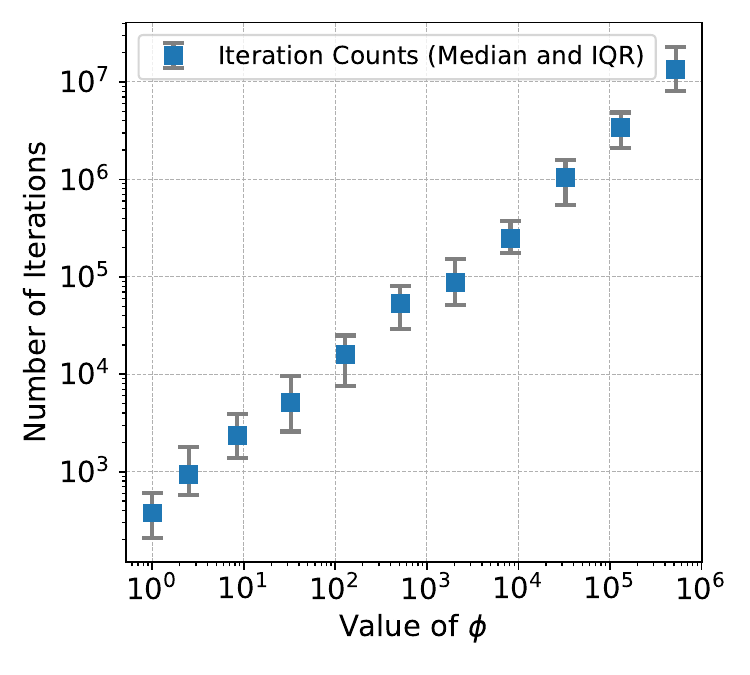}}
		\subcaptionbox{\small Stage II \vspace{-5pt}
		\label{fig:phi_stage2}}
			{\includegraphics[width=0.48\textwidth]{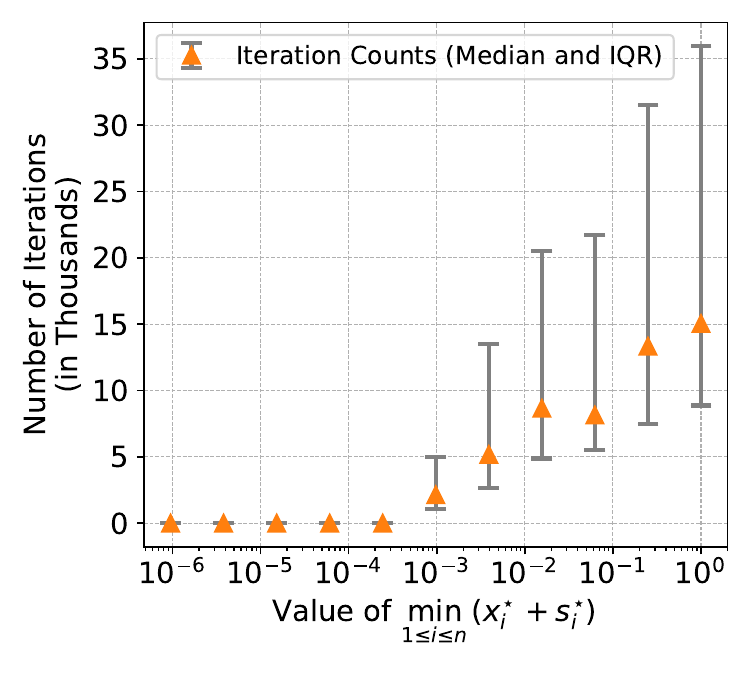}}
		\hspace*{\fill} % Add blank space on the right
		}
		{\normalsize The median values and the interquartile range (IQR) of Stage-I and Stage-II iteration counts for the random LP instances with various values of $\phi$ and $\min_{1\le i \le n} (x^\star_i + s^\star_i)$.\vspace{5pt}
		\label{fig:dependence_on_varphi}}
		{}
	\end{figure}

Now we confirm the difficulty of these instances via experiments. We follow the same experimental setup of rPDHG as in Section \ref{sec:experiments}, and we still consider an instance solved when a primal-dual solution pair is within Euclidean distance $10^{-4}$ of the optimal solution.  We set $m=50$ and $n = 100$. For each $l\in\{0,1,\dots,10\}$, we generate 100 LP instances and compute the first quartile, median, and third quartile of both Stage-I and Stage-II iteration counts of rPDHG.   
Figure \ref{fig:dependence_on_varphi} shows the relation between these statistics of the Stage-I iteration count and the value of $\phi$, and the relation between these statistics of the Stage-II iteration count and the value of $\min_{1 \le i \le n} (x_i^\star + s_i^\star)$, for each family of instances (grouped by $l$). 

Theorem \ref{thm:complexity-conditioned} predicts that the Stage-I iteration count grows linearly with the value of $\phi$.  If the theory is exact, in the log-log plot in Figure \ref{fig:phi_stage1}, the slopes of the three statistics should all be equal to $1$. 
Figure \ref{fig:phi_stage1} indeed shows that the three statistics of the Stage-I iteration count all have a slope of approximately $1$, confirming the linear dependence of the high-probability upper bound on $T_{basis}$ in Theorem \ref{thm:complexity-conditioned} on $\phi$.  This also validates our approach of controlling the value of $\frac{\frac{1}{n}\sum_{i=1}^{n}(\hat{x}_i + \hat{s}_i)}{\min_{1 \le i \le n} (\hat{x}_i + \hat{s}_i)}$ to generate challenging LP instances for rPDHG (and perhaps even other first-order LP methods). When $l = 10$, $\phi_l  \ge 5\times 10^5$ and the corresponding iteration count is likely to be higher than $10^7$. These instances are significantly more difficult than the LP instances studied in Section \ref{subsec:dependence_on_n_practice} for the high-probability performance of rPDHG. See the data points in Figure \ref{fig:dependence_on_n} along the horizontal axis of $n$ near $10^2$ for comparison.

Additionally, Theorem \ref{thm:complexity-conditioned} predicts that the Stage-II iteration count grows at most logarithmically with the value of $\min_{1 \le i \le n} (x_i^\star + s_i^\star)$ in a high-probability sense. This relationship can indeed be observed from Figure \ref{fig:phi_stage2}, particularly for $\min_{1 \le i \le n} (x_i^\star + s_i^\star) \ge 10^{-3}$.  
Theorem \ref{thm:complexity-conditioned} also indicates that $T_{local}$ is likely to be $0$ when $\min_{1 \le i \le n}(x^\star_i + s^\star_i)$ is smaller than or equal to the tolerance $\eps$. Indeed, Figure \ref{fig:phi_stage2} shows that when $\min_{1 \le i \le n} (x_i^\star + s_i^\star) \le 2 \cdot 10^{-4}$, all quartiles of the Stage-II iteration count become zero, suggesting that a sufficiently accurate solution is found even before the iterates settle on the optimal basis. However, in our instances, a smaller value of $\min_{1 \le i \le n} (x_i^\star + s_i^\star)$ does not indicate an easier problem, as it leads to larger values of $l$ and the corresponding $\phi_l$ (see \eqref{def:smallestnnz_l} and \eqref{def:phi_l}).

In the remainder of this section, we prove  Theorem \ref{thm:complexity-conditioned}.

\subsection{Proof of Theorem \ref{thm:complexity-conditioned}}\label{subsec:proof_of_subsec:effect_of_optimal_solution}
We first derive a conditional high-probability upper bound on $\kappa\Phi$.

\begin{lemma}\label{lm:kappaphi_unique}
	For the random LP defined in Definition \ref{def:random_lp}, for any $\delta \in (0,1]$, it holds that
	\begin{equation*}
		\Pr\left[ \kappa\Phi  \le \frac{c_3 n^{2.5}m^{0.5}\phi}{\delta(d+1)} \;\middle|\; E_{\mathrm{uniq}} \right] \ge  1- \delta - 4e^{-c_1 m}  - \left(\tfrac{1}{2}\right)^{d+1} \, .
	\end{equation*}
\end{lemma}
\begin{proof}{Proof.}
	For a fixed $\delta\in(0,1]$, let
	$
		G:=\left\{\kappa\|B^{-1}\|\|A\|\le \frac{c_3n^{1.5}m^{0.5}}{\delta(d+1)}\right\}\cap\{u>0\}.
	$
	Since $u>0$ almost surely, the tail bound \eqref{eqoflm:bound_of_kappaBinvA} gives
	\begin{equation*}
		\Pr[G]\ge 1-\delta-4e^{-c_1m}-\left(\tfrac{1}{2}\right)^{d+1}.
	\end{equation*}
	Recall the convention following \eqref{eq:def_varphi} that $\|B^{-1}\|:=\infty$ when $B$ is not full-rank. Thus, on $G$, $B$ is full-rank and $(\hat{x},\hat{s})$ satisfies strict complementary slackness. Lemma \ref{lm:full-rank-unique-optima} then implies that $G\subseteq E_{\mathrm{uniq}}$ and that $x^\star=\hat{x}$ and $s^\star=\hat{s}$ on $G$. Hence, by \eqref{eq:def_varphi} and \eqref{eq:def_of_phi}, $\varphi=n\phi$ on $G$. The definition of $G$ and Lemma \ref{lm:upperbound_kappa_phi_what} therefore yield
	\begin{equation*}
		\kappa\Phi\le \kappa\|B^{-1}\|\|A\|\varphi
		\le \frac{c_3n^{2.5}m^{0.5}\phi}{\delta(d+1)}
	\end{equation*}
	on $G$.
	Moreover, Lemmas \ref{lm:full-rank-unique-optima} and \ref{lm:full-rank-probability} imply $\Pr[E_{\mathrm{uniq}}]\ge 1-e^{-\nu m}>0$. Consequently,
	\begin{equation}  
		\Pr\left[\kappa\Phi\le \frac{c_3n^{2.5}m^{0.5}\phi}{\delta(d+1)}\;\middle|\;E_{\mathrm{uniq}}\right] \ge \frac{\Pr[G]}{\Pr[E_{\mathrm{uniq}}]}
		\ge \Pr[G]
		\ge 1-\delta-4e^{-c_1m}-\left(\tfrac{1}{2}\right)^{d+1}.
	\end{equation}
	where the second inequality uses $0 < \Pr[E_{\mathrm{uniq}}]\le1$.
	\hfill\Halmos\end{proof}

\begin{proof}{Proof of Theorem \ref{thm:complexity-conditioned}.}
	By \eqref{eq:basis_identification_T}, $T_{basis}\le\check{c}_1\kappa\Phi\ln(\kappa\Phi)$ for an absolute constant $\check{c}_1>1$. Furthermore, $\kappa\Phi\ge1$ by \eqref{def_geometric_measure} and $\kappa\ge1$. Since $z\mapsto z\ln z$ is increasing on $[1,\infty)$, Lemma \ref{lm:kappaphi_unique} proves \eqref{eq:optimal_basis_conditioned} with $C_0=c_1$ and $C_1=\check{c}_1c_3$.

	For the local convergence bound, Lemma \ref{lm:bound_of_BinvA} gives $\Pr[\|B^{-1}\|\|A\|\le c_2\sqrt{mn}/\delta]\ge1-\delta-2e^{-c_1m}$. This event implies that $B$ is full-rank, and its intersection with $\{u>0\}$ is therefore contained in $E_{\mathrm{uniq}}$. On this intersection, \eqref{eq:local_linear_convergence_T} gives the bound inside \eqref{eq:local_convergence_conditioned} with $C_2=\check{c}_2c_2$. Thus, the conditional probability in \eqref{eq:local_convergence_conditioned} is at least
	\begin{equation*}
		\frac{\Pr\left[\|B^{-1}\|\|A\|\le \tfrac{c_2\sqrt{mn}}{\delta},\ u>0\right]}{\Pr[E_{\mathrm{uniq}}]}
		\ge 1-\delta-2e^{-c_1m}.
	\end{equation*}
	Taking $C_0=c_1$, the constants $C_0,C_1,C_2$ depend only (and at most polynomially) on $\sigma_A$.
	\hfill\Halmos\end{proof}

% %\THEEndNotes
% \begingroup \parindent 0pt \parskip 4ex
% \def\enotesize{\normalsize} 
% \theendnotes
% \endgroup

% Appendix here
% Options are (1) APPENDIX (with or without general title) or
%             (2) APPENDICES (if it has more than one unrelated sections)
% Outcomment the appropriate case if necessary
%
% \begin{APPENDIX}{<Title of the Appendix>}
% \end{APPENDIX}
%
%   or
%
% \begin{APPENDICES}
% \section{<Title of Section A>}
% \section{<Title of Section B>}
% etc
% \end{APPENDICES}

% Acknowledgments here

\begin{APPENDICES}
    % Use distinct hyperlink destinations for appendix sections and their equations.
    \renewcommand*{\theHsection}{appendix.\Alph{section}}
    \SingleSpacedXI
    
	\section{Restarted PDHG Method for Linear Programming}\label{appsec:restartPDHG}
Algorithm \ref{alg: PDHG with restarts} presents the framework of rPDHG. 
\begin{algorithm}[htbp]
	\SetAlgoLined
	{\bf Input:} Initial iterate $(x^{0,0}, y^{0,0})=(0,0)$, $\ell \gets 0$, $k\gets 0 $, step sizes $\tau,\sigma$ \;
	\Repeat{\text{$(x^{\ell,0}, y^{\ell,0})$ satisfies some convergence condition }}{
	\Repeat{satisfying the restart condition}{
	\textbf{conduct one step of PDHG: }$(x^{\ell,k+1},y^{\ell,k+1}) \gets \textsc{OnePDHG}(x^{\ell,k},y^{\ell,k})$ \;\label{line:one-step-PDHG}\label{line:onepdhg}
	\textbf{compute the average iterate: }$(\bar{x}^{\ell,k+1},\bar{y}^{\ell,k+1})\gets\frac{1}{k+1} \sum_{i=1}^{k+1} (x^{\ell,i},y^{\ell,i})$
	\label{line:average} \;  \label{line:output-is-average-of-iterates}
	$k\gets k+1$ \;
	}\label{line:restart_condition}
	\textbf{restart the outer loop:} $(x^{\ell+1,0},y^{\ell+1,0})\gets (\bar{x}^{\ell,k},\bar{y}^{\ell,k})$, $\ell\gets \ell+1$, $k\gets 0$ \;
	}
	{\bf Output:} $(x^{\ell,0}, y^{\ell,0})$
	\caption{rPDHG: restarted-PDHG}\label{alg: PDHG with restarts}
\end{algorithm} 
It requires no matrix factorizations. Line \ref{line:onepdhg} of Algorithm \ref{alg: PDHG with restarts} is an iteration of the vanilla PDHG defined in \eqref{eq_alg: one PDHG}. Line \ref{line:output-is-average-of-iterates} can be computed efficiently by updating $(\bar{x}^{\ell,k},\bar{y}^{\ell,k})$ incrementally.
The restart condition in Line \ref{line:restart_condition} follows the $\beta$-restart criterion of \citet{applegate2023faster}, which triggers when the ``normalized duality gap'' of $(\bar{x}^{\ell,k},\bar{y}^{\ell,k})$ reduces to a $\beta$ fraction of $(x^{\ell,0},y^{\ell,0})$'s normalized duality gap. For the normalized duality gap's precise definition, see (4a) of \citet{applegate2023faster} for general primal-dual first-order methods and Definition 2.1 of \citet{xiong2024role} for PDHG in conic linear optimization.  The normalized duality gap measures the violations of feasibility and optimality, and $(x^{\ell,0},y^{\ell,0})$ becomes approximately optimal when its normalized duality gap is sufficiently small.

In practice, the restart condition is checked periodically (every few hundred iterations), and the normalized duality gap can be easily approximated by a separable norm variant within an absolute constant factor. See Section 6 of \citet{applegate2023faster} and Appendix A of \citet{xiong2024role} for its efficient approximation. This approach is implemented in \citet{applegate2023faster,lu2023cupdlpnew} and our experiments in Section \ref{sec:experiments}.

Following Theorem 3.1 of \citet{xiong2024accessiblenew}, we set the step sizes as:
$\tau = \frac{\lambda_{\min}}{2\lambda_{\max}}$ and $\sigma = \frac{1}{2\lambda_{\min}\lambda_{\max}}$,
where $\lambda_{\max}$ and $\lambda_{\min}$ denote the largest and smallest nonzero singular values of $A$, respectively. The restart parameter $\beta$ is set to $\frac{1}{e}$, where $e$ is the base of the natural logarithm. The largest singular value $\lambda_{\max}$ can be efficiently computed via power iterations. Using an imprecise estimate of $\lambda_{\min}$ is equivalent to applying rPDHG to a scalar-rescaled LP problem, so the method still converges to optimal solutions (\citet{applegate2023faster,xiong2024accessiblenew}). rPDHG with the above parameter setup has been studied in analyses of \citet{applegate2023faster,xiong2023computational,xiong2023relation,xiong2024role,xiong2024accessiblenew}.

\section{Proof of Lemma \ref{lm:full-rank-unique-optima}}\label{appsec:high-probability} 

\begin{proof}{Proof of Lemma \ref{lm:full-rank-unique-optima}.}

	First, observe that $\hat{x}$ and $\hat{s}$ are the optimal primal-dual solutions, so $\hat{x} \in \calX^\star$ and $\hat{s}\in \calS^\star$.

	For sufficiency, suppose $B$ is full-rank. Then $u^1\in\R^m$, $u^2\in \R^d$, $u^1 > 0$ and $u^2 > 0$  imply that $\hat{x}$ and $\hat{s}$ are both nondegenerate optimal solutions. Furthermore, since rows of $A$ are linearly independent, the instance has unique primal and dual optimal solutions: $\calX^\star = \{\hat{x}\}$ and $\calS^\star = \{\hat{s}\}$. Moreover, the corresponding $y$ satisfying $A^\top y + \hat{s} = c$ is also unique. Thus, $\calX^\star$, $\calY^\star$ and $\calS^\star$ are all singletons with $\calX^\star=\{\hat{x}\}$ and $\calS^\star = \{\hat{s}\}$.

	For necessity, we consider two cases. If $(\hat{x}, \hat{s})$ does not satisfy strict complementary slackness, at least a pair of strictly complementary optimal solutions exists. Alternatively, if $B$ is not full-rank but $u^1 > 0$ and $u^2 > 0$, then any solution $\tilde{x} = \left(\begin{smallmatrix}
        \tilde{x}_{[m]} \\ 0
    \end{smallmatrix}\right)$ satisfying $B \tilde{x}_{[m]} = B u^1$ and $\tilde{x}_{[m]} \geq 0$ is also optimal, where $[m]$ denotes $\{1,2,\dots,m\}$. Multiple such $\tilde{x}_{[m]}$ exist because $B$ is not full-rank and $u^1 > 0$. Therefore, in either case, the optimal solution sets $\calX^\star$ and $\calS^\star$ cannot both be singletons with $\calX^\star=\{\hat{x}\}$ and $\calS^\star = \{\hat{s}\}$.
\hfill\Halmos\end{proof}

\section{Proof of Section \ref{subsec:helper_lemmas}}\label{appsec:prove_thm_random}

\subsection{Proof of Lemma \ref{lm:max_singular_value_rectangular}}

We actually have the following more general result.
\begin{lemma}\label{lm:max_singular_value_rectangular_base}
	Let $A$ be a random matrix as defined in Definition \ref{def:random matrix}, where $A \in \R^{m \times n}$ and $n \ge m$. Then, for all $t \ge 10\sigma_A$, we have:
	\begin{equation}\label{eqoflm:max_singular_value_rectangular}
		\Pr\left[\sigma_{1}(A) \geq t\sqrt{n} \right] \le \exp\left(-\tfrac{t^2 n}{16\sigma_A^2}\right)  \ .
	\end{equation} 
\end{lemma}

\begin{proof}{Proof.}
We apply the argument of Proposition 2.3 of
\citet{rudelson2009smallest} to $A^\top$, so their dimensions
$(N,n)$ correspond to our $(n,m)$. First, replace both $1/2$-nets
in their proof by $1/4$-nets. The volume comparison in the proof
of their Proposition 2.1 gives cardinality bounds $9^n$ and $9^m$,
respectively.

Next, by applying Lemma \ref{lm:hoeffding} to both tails, the
constants $c_1$ and $C_1$ in the second displayed equation of
their proof can be taken as $1/(2\sigma_A^2)$ and $2$,
respectively, because the squared coefficients in the bilinear
form sum to one for any fixed pair of unit vectors.

Finally, with the modified nets, Proposition 2.2 therein gives
an approximation factor $(1-1/4)^{-2}=16/9<2$. Thus, in the
final union-bound step, we apply the fixed-pair tail bound at
threshold $t\sqrt{n}/2$, rather than $t\sqrt{n}$. This replaces
the last displayed estimate of their proof by
\begin{equation}
\Pr\left[\sigma_1(A)\ge t\sqrt{n}\right]
\le 2\cdot 9^n\cdot 9^m
\exp\left(-\frac{t^2n}{8\sigma_A^2}\right)
\le \exp\left(-\frac{t^2n}{16\sigma_A^2}\right),
\end{equation}
where the last inequality follows from
$
\ln 2+(m+n)\ln 9
\le n\ln 162
<\frac{25n}{4}
\le\frac{t^2n}{16\sigma_A^2}
$
for $m\le n$ and $t\ge 10\sigma_A$. This proves the lemma.
\hfill\Halmos
\end{proof} 

Now we can prove Lemma \ref{lm:max_singular_value_rectangular}.

\begin{proof}{Proof of Lemma \ref{lm:max_singular_value_rectangular}.}
	Setting $t=10\sigma_A$ in \eqref{eqoflm:max_singular_value_rectangular} gives
	$\Pr\left[\sigma_1(A)\ge10\sigma_A\sqrt{n}\right]\le e^{-25n/4}\le e^{-6n}$.
\hfill\Halmos\end{proof}

\subsection{Proof of Lemma \ref{lm:tail_of_product}} 
We first present a technical result before proving Lemma \ref{lm:tail_of_product}.

\begin{lemma}\label{lm:bound_of_deltalndelta}
	For any $y \in (0,\tfrac{1}{e}]$, if there exists $x \in (0,\tfrac{1}{e}]$ such that $x\ln(\tfrac{1}{x}) = y$, then $x \ge \frac{y}{2\ln(1/y)}$.
\end{lemma}

\begin{proof}{Proof.}
	Define $f(x):=x\cdot \ln(1/x)$ and $x_0 := \frac{y}{2\ln(1/y)}$.  Then we have $f(x_0) = \frac{y}{2\ln(1/y)}\cdot \ln(2\ln(1/y)/y) = \frac{y}{2\ln(1/y)}\cdot [\ln(2\ln(1/y)) + \ln(1/y)]= \frac{y}{2} + \frac{y}{2}\cdot \frac{\ln(2\ln(1/y))}{\ln(1/y)}$. Note that as $y\in(0,1/e]$, $2\ln(1/y) - 1/y\le \max_{t\ge e} \left[2\ln(t) -t\right]= 2-e < 0$, so $\frac{\ln(2\ln(1/y))}{\ln(1/y)} \le 1$ and thus $f(x_0) \le \frac{y}{2} + \frac{y}{2} = y$.
	Finally, since $f(x):=x\cdot \ln(1/x)$ is monotonically increasing on $(0,1/e]$, we have $x_0 \le x$. This proves the statement.
\hfill\Halmos\end{proof}

Now we prove Lemma \ref{lm:tail_of_product}.
\begin{proof}{Proof of Lemma \ref{lm:tail_of_product}.}
	First of all, we prove that for all $T > C_1C_2$:
	\begin{equation}\label{eq:tail_of_product_2}
		\Pr[XY \ge T] \le \tfrac{2C_1C_2}{T} + \tfrac{C_1C_2}{T}\cdot \ln\left(\tfrac{T}{C_1C_2}\right) + \delta_1 + 2 \delta_2 \ .
	\end{equation}
	For a random variable $\eta$, let $F_\eta$ denote the distribution function of $\eta$.
	The probability $\Pr[XY\ge T]$ has the following upper bound:
	\begin{equation}\label{eq:tail_of_product_7}
		\begin{aligned}
			\Pr[XY\ge T] = \int_0^\infty \Pr\left[X\ge \tfrac{T}{y}\right]\mathrm{d} F_Y(y) \overset{\eqref{eq:tail_of_product_1}}{\le}  \delta_1 + \underbrace{\int_0^\infty  \min\left\{\tfrac{C_1y}{T},1\right\} \cdot \mathrm{d} F_Y(y)}_{\text{Denoted by $\mathbf{I}$}}    \ ,
		\end{aligned}
	\end{equation}
	in which
	\begin{equation}\label{eq:tail_of_product_8}
		\begin{aligned}
			\mathbf{I}     = \underbrace{\int_0^{\tfrac{T}{C_1}} \tfrac{C_1y}{T} \cdot \mathrm{d} F_Y(y)}_\text{Denoted by $\mathbf{II}$}  + \int_{\tfrac{T}{C_1}}^\infty \mathrm{d} F_Y(y)   =\mathbf{II} + \Pr\left[Y \ge \tfrac{T}{C_1}\right] \overset{\eqref{eq:tail_of_product_1}}{\le} \mathbf{II} + \tfrac{C_1C_2}{T} + \delta_2 \ .
		\end{aligned}
	\end{equation}
	The term $\mathbf{II}$ is $\tfrac{C_1}{T}\int_0^{\tfrac{T}{C_1}} y \cdot \mathrm{d} F_Y(y)$, and $\tfrac{T}{C_1}\cdot \mathbf{II}$ has the following upper bound:
	\begin{equation}\label{eq:tail_of_product_9}{\small
			\begin{aligned}
				\tfrac{T}{C_1}\cdot \mathbf{II} & = \int_0^{\tfrac{T}{C_1}} y \cdot \mathrm{d} F_Y(y) =  y\cdot F_Y(y)\bigg|_0^{\tfrac{T}{C_1}} - \int_0^{\tfrac{T}{C_1}}  F_Y(y) \cdot \mathrm{d} y                                                                                 \\
				                                & = \tfrac{T}{C_1}\cdot \Pr\left[Y \le \tfrac{T}{C_1}\right] - \int_0^{\tfrac{T}{C_1}}  \left(1-\Pr\left[Y > y\right]\right)  \mathrm{d} y  \le \tfrac{T}{C_1}  - \tfrac{T}{C_1} + \int_0^{\tfrac{T}{C_1}}   \Pr[Y > y]  \mathrm{d} y \\
				                                &
				\overset{\eqref{eq:tail_of_product_1}}{\le}    \int_0^{\tfrac{T}{C_1}}   \min\left\{\tfrac{C_2}{y}+\delta_2,1\right\}  \mathrm{d} y  \le   \tfrac{T}{C_1}\cdot \delta_2 + \int_0^{\tfrac{T}{C_1}}   \min\left\{\tfrac{C_2}{y},1\right\}  \mathrm{d} y                \\
				                                & \le \tfrac{T}{C_1}\cdot \delta_2 +  \int_0^{C_2}   1 \mathrm{d} y  +
				\int_{C_2}^{\tfrac{T}{C_1}}    \tfrac{C_2}{y}  \mathrm{d} y  = \tfrac{T}{C_1}\cdot \delta_2 + C_2 + C_2\ln\left(\tfrac{T}{C_1C_2}\right) \ .
			\end{aligned}}
	\end{equation}
	In other words,
	\begin{equation}\label{eq:tail_of_product_11}
		\mathbf{II}\le   \tfrac{C_1}{T}\left(\tfrac{T}{C_1}\cdot \delta_2+ C_2 + C_2\ln\left(\tfrac{T}{C_1C_2}\right)\right)
		=  \delta_2 + \tfrac{C_1C_2}{T} + \tfrac{C_1C_2}{T}\cdot \ln\left(\tfrac{T}{C_1C_2}\right)
		\ .
	\end{equation}
	Finally, substituting the upper bound of $\mathbf{II}$ in \eqref{eq:tail_of_product_11} back into \eqref{eq:tail_of_product_8} and \eqref{eq:tail_of_product_7} yields \eqref{eq:tail_of_product_2}.

	We are particularly interested in \eqref{eq:tail_of_product_2} when $T \ge 3C_1 C_2$, as otherwise the bound is trivial. We use $\delta_0$ to denote $\tfrac{C_1C_2}{T}$ (which takes values in $(0,1/3]$), and \eqref{eq:tail_of_product_2} reduces to $\Pr\left[XY \ge \tfrac{C_1C_2}{\delta_0}\right] \le 3\delta_0 \cdot \ln(\tfrac{1}{\delta_0}) + \delta_1 + 2 \delta_2$.  We let $\hat{\delta}:=\delta_0\cdot \ln(1/\delta_0)$, and $\hat{\delta}$ may take any value in $(0,\frac{\ln(3)}{3}]$. Using Lemma \ref{lm:bound_of_deltalndelta}, $\delta_0$ is at least as large as $\frac{\hat{\delta}}{2\ln(1/\hat{\delta})}$. Therefore,
	\begin{equation}\label{eq:tail_of_product_12}
		\Pr\left[XY \ge \tfrac{2C_1C_2 \cdot \ln(1/\hat{\delta})}{\hat{\delta}}\right]\le \Pr\left[XY \ge \tfrac{C_1C_2}{\delta_0}\right]\le 3\hat{\delta} + \delta_1 + 2 \delta_2 \ .
	\end{equation}
	Finally, we replace $3\hat{\delta}$ by  $\delta$ (which may take any value in $(0,1]$) and \eqref{eq:tail_of_product_12} then becomes the desired inequality \eqref{eq:tail_of_product_1_2}.
\hfill\Halmos\end{proof}

\section{Proof of Section \ref{subsec:useful_helper_lemmas_gaussian}}\label{appsec:useful_helper_lemmas_gaussian}

\subsection{Proof of Lemma \ref{lm:distribution_of_sigmai}}

\begin{proof}{Proof of Lemma \ref{lm:distribution_of_sigmai}.}
	According to the min-max principle for singular values, the $k$-th largest singular value $\sigma_k(W)$ has the following equivalent representation:
	\begin{equation}\label{eq:distribution_of_sigmai_1}
		\sigma_k(W)  = \max_{S:\dim(S) = k}  \min_{x\in S, \|x\| = 1} \|Wx\| \ .
	\end{equation}
	Let $S_k := \{x\in\R^n: x_{k+1}=x_{k+2}=\cdots=x_n=0\}$ denote the $k$-dimensional subspace of vectors with zeros in the last $n-k$ components. Since $\dim(S_k) = k$, \eqref{eq:distribution_of_sigmai_1} implies that $\sigma_k(W)$ has the following lower bound:
	\begin{equation}\label{eq:distribution_of_sigmai_2}
		\sigma_k(W)  \ge  \min_{x\in S_k, \|x\| = 1} \|Wx\| = \min_{y\in\R^{k},\|y\|= 1}  \|W_{[k]}y\|  = \sigma_k(W_{[k]}) 
	\end{equation}
	where $W_{[k]}$ denotes the first $k$ columns of $W$. Thus, $\sigma_k(W)$ is bounded below by the smallest singular value of $W_{[k]}$.

	Since $W_{[k]}$ is also a sub-Gaussian matrix, it further holds that 
	\begin{equation}\label{eq:distribution_of_sigmai_3}
		(C_{rv1} \delta_0)^{n-k+1} + e^{-C_{rv2} n} \ge \Pr\left[\sigma_k(W_{[k]}) < \delta_0(\sqrt{n}-\sqrt{k-1})\right] \ge  \Pr\left[\sigma_k(W_{[k]}) < \delta_0\cdot \frac{n-k+1}{2\sqrt{n}}\right]
	\end{equation}
	for all $\delta_0 > 0$, where the first inequality follows from Lemma \ref{lm:min_singular_value_rectangular} applied to $W_{[k]}$ and the second inequality is due to $\sqrt{n} - \sqrt{k-1} = \frac{n-k+1}{\sqrt{n} + \sqrt{k-1}} \ge \frac{n-k+1}{2\sqrt{n}}$. Therefore,  this yields the following bound on all singular values for all $\delta_0 > 0$:
	\begin{equation}\label{eq:distribution_of_sigmai_5}
		\begin{aligned}
			& \Pr\left[\sigma_k(W) \ge \frac{\delta_0 (n-k+1)}{2\sqrt{n}}\text{ for }k=1,2,\dots,n\right]  \\
			 \ge \ & 1 - \sum_{k=1}^n \Pr\left[\sigma_k(W) < \frac{\delta_0 (n-k+1)}{2\sqrt{n}}\right] 
			\overset{\eqref{eq:distribution_of_sigmai_2}}{\ge} 1 - \sum_{k=1}^n \Pr\left[\sigma_k(W_{[k]}) < \frac{\delta_0 (n-k+1)}{2\sqrt{n}}\right]  \\
			\overset{\eqref{eq:distribution_of_sigmai_3}}{\ge} & 1-   \sum_{k=1}^{n}(C_{rv1} \delta_0)^{n-k+1} - n\cdot e^{-C_{rv2} n}  \ge 1 - 2C_{rv1}\delta_0 -n \cdot e^{-C_{rv2} n} \ .
		\end{aligned}
	\end{equation}
	The final inequality holds because for $\delta_0 \in (0,\frac{1}{2C_{rv1}})$, $\sum_{k=1}^{n}(C_{rv1} \delta_0)^k \le \sum_{k=1}^{\infty}(C_{rv1} \delta_0)^k = \frac{C_{rv1} \delta_0}{1-C_{rv1} \delta_0} \le 2C_{rv1} \delta_0$. Replacing $2C_{rv1} \delta_0$ in \eqref{eq:distribution_of_sigmai_5} with $\delta$ establishes \eqref{eq_of_lm:distribution_of_sigmai} for all $\delta \in (0,1)$. As for the case $\delta \ge 1$, \eqref{eq_of_lm:distribution_of_sigmai} is trivial.

\hfill\Halmos\end{proof}

\subsection{Proof of Lemma \ref{lm:tail_of_Binvv}}

\begin{proof}{Proof of Lemma \ref{lm:tail_of_Binvv}.}
	Let $M$ denote $W^{-\top}W^{-1}$. Then $\|W^{-1}v\|^2 = v^\top W^{-\top}W^{-1}v = v^\top M v$.
	We will first prove the following tail bound of $v^\top M v$ for all $\gamma \ge 0$:
	\begin{equation}\label{eq:lm:tail_of_Binvv_7}
\begin{aligned}
& \Pr\left[v^\top M v \ge \E\left[v^\top M v\right] + \frac{2\gamma\sigma^2 n }{c_0^2}\right]
\le 2\cdot \exp\!\left(-2C_{hw}\min\{\gamma^2,\gamma\}\right) \ .
\end{aligned}
\end{equation}
	
	To establish \eqref{eq:lm:tail_of_Binvv_7}, we apply Lemma \ref{lm:Hanson-Wright}, which requires bounds on $\|M\|$ and $\|M\|_F$. We begin by characterizing the singular values of $M$. For $k=1,2,\dots,n$:
	\begin{equation}\label{eq:lm:tail_of_Binvv_1}
		\sigma_k(M) = \sigma_k(W^{-\top}W^{-1}) =  \frac{1}{\sigma_{n+1-k}^2(W)}
	\end{equation}
	Then we have the following upper bounds for $\|M\|_F^2$ and $\|M\|$:
	\begin{equation}\label{eq:lm:tail_of_Binvv_2}
		\|M\| = \sigma_1(M) \overset{\eqref{eq:lm:tail_of_Binvv_1}}{=}  \frac{1}{\sigma_n^2(W)} \overset{\eqref{eqoflm:tail_of_Binvv_1}}{\le} \frac{n}{c_0^2}
 	\end{equation}
	\begin{equation}\label{eq:lm:tail_of_Binvv_3}
		\|M\|_F^2 = \sum_{i=1}^n\sigma_i^2(M)  \overset{\eqref{eq:lm:tail_of_Binvv_1}}{=}  \sum_{i=1}^n  \frac{1}{\sigma_i^4(W)} \overset{\eqref{eqoflm:tail_of_Binvv_1}}{\le} 
		\sum_{i=1}^n \frac{n^2}{c_0^4 i^4} \le \frac{n^2}{c_0^4}\sum_{i=1}^\infty \frac{1}{i^4} = \frac{n^2}{c_0^4} \cdot \frac{\pi^4}{90} \le \frac{2n^2}{c_0^4} \ .
	\end{equation}
	where the last equality uses the definition of the Riemann zeta function and its value at $4$. 
	With the above upper bounds of $\|M\|_F^2$ and $\|M\|$, now we can use Lemma \ref{lm:Hanson-Wright}:
	\begin{equation}\label{eq:lm:tail_of_Binvv_6}
		\begin{aligned}
		  \Pr\Big[v^\top M v \ge \E\left[v^\top M v\right] + t\Big]  \le 2\cdot \exp \left[-C_{hw} \cdot \min \left(\frac{c_0^4t^2}{2\sigma^4 n^2}, \frac{c_0^2t}{\sigma^2 n}\right)\right] \ .
		\end{aligned}
	\end{equation} 
	For any $\gamma \ge 0$, set $t = \frac{2\gamma\sigma^2 n}{c_0^2}$. Substituting this choice into \eqref{eq:lm:tail_of_Binvv_6} yields:
\begin{equation}\label{eq:lm:tail_of_Binvv_6_0}
\begin{aligned}
\Pr\left[v^\top M v \ge \E\left[v^\top M v\right] + \frac{2\gamma\sigma^2 n}{c_0^2}\right]
&\le 2\cdot \exp \left[-C_{hw} \cdot \min \left(\frac{c_0^4}{2\sigma^4 n^2}\cdot \frac{4\gamma^2\sigma^4 n^2}{c_0^4},
\frac{c_0^2}{\sigma^2 n}\cdot \frac{2\gamma\sigma^2 n}{c_0^2}\right)\right] \\
&= 2\cdot \exp \left[-2C_{hw} \cdot \min \left\{\gamma^2,\gamma\right\}\right] \ .
\end{aligned}
\end{equation}
This proves \eqref{eq:lm:tail_of_Binvv_7}.

	We now use \eqref{eq:lm:tail_of_Binvv_7} to derive the tail bound of $v^\top M v$. We start from the following expression of $\E\left[v^\top M v\right]$:
	\begin{equation}\label{eq:lm:tail_of_Binvv_4}
		\E\left[v^\top M v\right]  =  \E[\|W^{-1}v\|^2]  
		 =\sum_{i=1}^n \E\left[(W^{-1}v)_i^2\right]=\sum_{i=1}^n\left\|(W^{-1})_{i,\cdot}\right\|^2=\|W^{-1}\|_{F}^2 = \sum_{i=1}^n \frac{1}{\sigma_i^2(W)} 
	\end{equation}  
	where the third equality is because components of $v$ are independent and have zero mean and unit variance. 
	Furthermore, $\E\left[v^\top M v\right]$ is upper bounded as follows:
	\begin{equation}\label{eq:lm:tail_of_Binvv_5}
		\E\left[v^\top M v\right] \overset{\eqref{eq:lm:tail_of_Binvv_4}}{=} \sum_{i=1}^{n} \frac{1}{\sigma_i^2(W)} \overset{\eqref{eqoflm:tail_of_Binvv_1}}{\le} \sum_{i=1}^{n} \frac{n}{c_0^2 i^2} \le \frac{n}{c_0^2} \sum_{i=1}^\infty \frac{1}{i^2} = \frac{n}{c_0^2} \cdot \frac{\pi^2}{6}   \le    \frac{2n}{c_0^2} \ ,     
	\end{equation}  
	where the last equality uses the definition of the Riemann zeta function and its value at $2$. Therefore, since $\|W^{-1}v\|^2 = v^\top M v$, for all $t > 0$:
	\begin{equation}
		\Pr\left[\|W^{-1}v\|^2 \ge \frac{2n}{c_0^2}  + t\right] \le \Pr\left[v^\top M v \ge \E\left[v^\top M v\right]  + t\right] \ .
	\end{equation}
	Finally, letting $t = \frac{2\gamma \sigma^2 n}{c_0^2}$ and applying \eqref{eq:lm:tail_of_Binvv_7} to this inequality completes the proof.
\hfill\Halmos\end{proof}

\end{APPENDICES}

    \ACKNOWLEDGMENT{The author thanks Yinyu Ye, Moïse Blanchard and Haihao Lu for the helpful discussions.
	The author is also grateful to Robert M. Freund for the valuable suggestions while preparing this manuscript. The author also thanks the three anonymous referees for helpful comments and suggestions that improved the quality of the manuscript. Part of this work was performed at the Massachusetts Institute of Technology and was supported by AFOSR Grant No. FA9550-22-1-0356. Another part of this work was performed at Georgia Institute of Technology and was supported by ONR Award \#N00014-25-1-2088.}

% References here (outcomment the appropriate case)

% CASE 1: BiBTeX used to constantly update the references
%   (while the paper is being written).
%\bibliographystyle{informs2014} % outcomment this and next line in Case 1
%\bibliography{<your bib file(s)>} % if more than one, comma separated

\bibliographystyle{informs2014} % outcomment this and next line in Case 1
\bibliography{referencenew} % if more than one, comma separated

% CASE 2: BiBTeX used to generate mypaper.bbl (to be further fine tuned)
%\input{mypaper.bbl} % outcomment this line in Case 2

%If you don't use BiBTex, you can manually itemize references as shown below.

%\bibliographystyle{nonumber}

% \begin{thebibliography}{3}
% \providecommand{\natexlab}[1]{#1}
% \providecommand{\url}[1]{\texttt{#1}}
% \providecommand{\urlprefix}{URL }

% \bibitem[{Smith(2005)}]{smith2005}
% Smith J (2005) Optimal resource allocation in humanitarian logistics.
%   \emph{Journal of Operations Research} 30(2):123--135.
  
% \bibitem[{Jones(2010)}]{jones2010}
% Jones S (2010) Stochastic programming models for humanitarian logistics.
%   \emph{INFORMS Mathematics of Operations Research} 35(4):567--580.

% \bibitem[{Brown(2015)}]{brown2015}
% Brown D (2015) \emph{Introduction to Stochastic Programming} (Springer).

% \end{thebibliography}

%%%%%%%%%%%%%%%%%
\end{document}